\documentclass[11 pt]{amsart}

\usepackage{amscd, amssymb}

\newtheorem{pr}{Proposition}
\newtheorem{lm}{Lemma}

\newcommand{\proj}{\mathbf P}

\newcommand{\rarr}{\rightarrow}
\newcommand{\oh}{{\mathcal{O}}}
\newcommand{\com}{\mathbb{C}}
\newcommand{\Q}{\mathbb{Q}}
\newcommand{\Z}{\mathbb{Z}}
\newcommand{\R}{\mathbf{R}}

\newcommand{\lan}{\langle}
\newcommand{\ran}{\rangle}
\newcommand{\eqq}{\stackrel{\sim}{=}}
\newcommand{\deli}{\bigtriangleup}

\newcommand{\grad}{\nabla}

\newcommand{\gradd}{\grad_{\hb}} 
\newcommand{\graddi}{\grad_{\hb,i}}
\newcommand{\hb}{\hbar}

\def\scup{\mathbin{\text{\scriptsize$\cup$}}}

\newcommand{\bpf}{\noindent {\em Proof.} }
\newcommand{\epf}{\qed \vspace{+10pt}}

\begin{document}
\title{Rational curves on hypersurfaces \\
(after A. Givental)}
\author{R. Pandharipande}
\date{ 21 June 1998}
\maketitle

\pagestyle{plain}
\setcounter{section}{-1}
\section{\bf{Introduction}}
We describe here a remarkable relationship studied by
Givental between hypergeometric series and the quantum
cohomology of hypersurfaces in projective space [G1].
As the quantum product involves genus 0 Gromov-Witten
invariants, a connection between hypergeometric series and the
geometry of rational curves on the hypersurfaces is made.  
While the most general context for such relationships
has not yet been understood, 
analogous results for complete
intersections in smooth toric varieties and flag varieties 
have been pursued by Givental [G2] and Kim [Ki].

The first discovery in this subject was the
startling prediction from Mirror symmetry 
by Candelas, de la Ossa, Green, and
Parkes [COGP] of the
numbers of rational curves on quintic $3$-folds in $\proj^4$.
We recount an equivalent form of their original prediction. 
Let $I_i(t)$ be defined by:
$$\sum_{i=0}^3 I_i H^i = \sum_{d=0}^{\infty} e^{(H+d)t} 
\frac{\Pi_{r=1}^{5d} (5H+r)}{\Pi_{r=1}^d(H+r)^5} \ \ \text{mod} \ H^4.$$
The functions $I_i(t)$ are a basis of solutions of the
Picard-Fuchs differential equation
$$\Big( \frac{d}{dt}\Big)^4 I= 5 e^t \Big(5 \frac{d}{dt}+1\Big)
\Big(5 \frac{d}{dt}+2\Big) \Big(5 \frac{d}{dt}+3\Big) 
\Big(5 \frac{d}{dt}+4\Big) I$$
arising in the $B$-model 
from the variation of Hodge structures of a specific
family of Calabi-Yau 3-folds.
Let $n_d$ be the virtual number of degree $d$
rational curves on a general quintic $3$-fold in $\proj^4$.
Let the change of variables 
$T(t)= I_1/I_0 \ (t)$ define a new coordinate $T$. 
The functions  $J_i = I_i/I_0 \ (T)$ in the new variable 
were predicted to equal an $A$-model series:
\begin{equation}
\label{msy}
\sum_{i=0}^3 J_i  H^i =
e^{HT} + \frac{H^2}{5} \sum_{d=1}^{\infty} n_d d^3
\sum_{k=1}^{\infty} \frac{e^{(H+kd)T}}{(H+kd)^2} \ \ \text{mod} \ H^4
\end{equation}
and satisfy the differential equation:
\begin{equation}
\label{ddeq}
\frac{d^2}{dT^2} \frac{1}{K(e^T)} \frac{d^2}{dT^2} \ J_i =0,
 \ \ \ \text{where} \ \ \  
K(e^T)= 5 + \sum_{d=1}^{\infty} n_d d^3 \frac{e^{dT}}
{1-e^{dT}}.
\end{equation}
As the enumerative
geometry of quintic 3-folds was not known  
to have any structure at all, these formulas were completely
unexpected. 
There is a large literature in
both physics and mathematics on 
Mirror symmetry for Calabi-Yau 3-folds. The $A$-model / $B$-model
framework is described in [W1]. A mathematical perspective
can be found in [Mo], [CK]. Higher dimensional Calabi-Yau
manifolds are considered in [GMP].

The numbers $n_d$ in (\ref{msy}) and (\ref{ddeq})
have the following interpretation:
if the rational curves on the
general quintic $3$-fold $Q$ were 
nonsingular, isolated with
balanced normal bundle, and disjoint, then
$n_d$ would 
simply be the number of degree $d$ rational curves on $Q$. 
However, this strong assumption is false [V] -- there exist
nodal degree 5 rational curves on $Q$. The nonexistence of
families of rational curves in $Q$ is still open (Clemens' conjecture).
A mathematically precise statement of conjecture (\ref{msy}) 
requires
a substantial program to define $n_d$: moduli spaces of maps, their
virtual classes, and Gromov-Witten invariants. This program 
has been completed in both symplectic and algebraic geometry
through the recent work of many mathematicians (see
the foundational papers [RT], [KM]). 
The virtually enumerative numbers $n_d$ are defined via the 
genus 0 Gromov-Witten invariants $N_d$ of the quintic 
by a formula accounting for 
multiple cover contributions [AM], [M]:
\begin{equation}
\label{manmul}
\sum_{d=1}^{\infty} N_d q^d = \sum_{d=1}^{\infty} \sum_{k=1}
^{\infty} {n_d}{k^{-3}} q^{kd}.
\end{equation}
The numbers $n_d$  are enumerative at least for
$d\leq 9$
 [K], [KJ1]. An outlook in
higher degrees may be found in [KJ2].

The central relationship in Givental's work may be explained as follows.
Let $X$ be a hypersurface in $\proj^m$ of degree
$l \leq m+1$. A correlator $S_X$ is defined
via the quantum product and related quantum 
differential equations associated
to $X$ (see Section 2). 
$S_X$ is closely related to the 
hypergeometric series:
\begin{equation}
\label{hgs}
S^*_X=\sum_{d=0}^{\infty} e^{(H+d)t} 
\frac{\Pi_{r=1}^{ld} (lH+r)}{\Pi_{r=1}^d(H+r)^{m+1}} \ \text{mod} \ H^{m}.
\end{equation}
The precise relationship is divided in 3 cases.
\begin{enumerate}
\item[(i)]
If $l<m$,  then $S^*_X= S_X$. 
\item[(ii)] If $l=m$, then $e^{-m!e^t} S^*_X= S_X$.
\item[(iii)] If $l=m+1$, then $S^*_X$ and $S_X$ are related
              by an explicit transformation (see Section 4).
\end{enumerate}
In the case of the quintic $3$-fold, $S_X$ is exactly the right
side of equation (\ref{msy}) (see Section 4.5). 
The transformation (iii) then specializes to the
Mirror symmetry prediction proving (\ref{msy}).
Equation (\ref{ddeq}) is a consequence of the
quantum differential equation. 
The results in cases (i) and (ii) have direct applications to
the quantum cohomology ring of the corresponding hypersurfaces
(see Section 3).
Givental has
suggested that cases (i) and (ii) correspond to
non-compact Mirrors.

The plan of the paper is as follows. Section 1 contains
a rapid review of Gromov-Witten invariants,
descendents,  quantum products, and 
quantum differential equations. 
In Section 2, a new quantum product $*_X$ is defined
on the cohomology of the ambient space $\proj^m$.
The $*_X$-product
greatly clarifies the relationship between quantum structures
on $X$ and $\proj^m$.
Givental's correlator $S_X$ is naturally defined
via differential equations arising from the $*_X$-product. 
This product appears in [G1], [Ki] and was explained to the
author by T. Graber.
Section 3 covers cases (i) and (ii) where
$l\leq m$. These are much easier than the Calabi-Yau
case which is established in Section 4.
The treatment in Sections 3 and 4
follows [G1] with some augmentation and 
modification.

The main technical tool needed in Givental's 
approach is an explicit 
localization formula in equivariant cohomology for the
natural torus action on the moduli space of maps 
$\overline{M}_{0,n}(\proj^m, d)$. As this moduli space
is a nonsingular stack, the Bott residue formula holds.
The fixed point loci of the
torus action as well as the precise equivariant normal bundle
determinations have been explained in detail in [Ko].
For Givental's arguments in the smooth toric case [G2], 
a virtual localization formula [GP] is necessary as the
moduli space of maps may be quite ill-behaved.  

The number $n_1=2875$ of lines on a general quintic $3$-fold
was obtained in the $19^{th}$ century by Schubert (via 
intersection calculations 
in the Grassmannian $\mathbf{G}(\proj^1,
\proj^4)$).
The numbers $n_2= 609250$ and
$n_3=317206375$ of conics and twisted cubics  were computed by S. Katz [K]
and Ellingsrud and Str\o mme [ES1] respectively.  
Localization was first applied to the enumerative geometry of
quintics  in [ES2].
The method of torus localization on $\overline{M}_{0,n}(\proj^m,d)$
was developed by Kontsevich in [Ko] precisely to attack the
Mirror prediction. The resulting formulas determined all the numbers
$n_d$ by a complex sum over graphs. This summation yielded the
first mathematical computation of $n_4=242467530000$ [Ko]. 
An important
aspect of the 
argument in Sections 3 and 4 is an organization of graph sums.

A complete proof of the Mirror prediction for
quintics by Lian, Liu, and Yau
using localization formulas 
has appeared recently in [LLY]. 
The argument announced by Givental in [G1] 
yields a complete proof of (i)-(iii). 
It is the
latter proof that is explained here (see [LLY], [G3] for
a comparison of viewpoints). Givental's work is also discussed in
[CK] and [BDPP].
So far, mathematical 
approaches to the $A$-model series 
do not involve the $B$-model at all.
While these results verify the predictions
of [COGP], the full correspondence of
Mirror symmetry remains to be mathematically explained.

The author wishes to thank G. Bini, A. Elezi,  
W. Fulton, E. Getzler,
T. Graber, S. Katz, S. Kleiman, 
B. Kim, M. Polito, and M. Thaddeus for many conversations on Givental's work.
The author was partially supported by a National Science
Foundation post-doctoral fellowship. Thanks are also
due to the Mittag-Leffler Institute where the author
learned much of this material. This paper was completed
at the Scuola Normale Superiore di Pisa.

\section{\bf {The quantum differential equation}}
\subsection{Descendents}
A nonsingular algebraic variety $X$ is {\em convex} if
$$H^1(\proj^1, \mu^*T_X)=0$$ for all maps
$\mu: \proj^1 \rarr X$. The main examples
of compact convex varieties are $X=\mathbf{G}/\mathbf{P}$ where 
$\mathbf{G}$ is a linear algebraic group and $\mathbf{P}$ is a
parabolic subgroup. The case of most interest here
is $X=\proj^m$. 
The space 
$\overline{M}_{0,n}(X, \beta)$
of
$n$-pointed genus 0  stable maps representing
the class $\beta\in H_2(X, \mathbb{Z})$ is 
a coarse moduli space with  quotient
singularities (or a nonsingular Deligne-Mumford stack)
of pure dimension $$\text{dim}(X)+\int_\beta c_1(T_X)
+n-3$$
in case $X$ is convex [Ko], [FP].
Let $e_i: \overline{M}_{0,n}(X, \beta)
\rarr X$ be the $i^{th}$ evaluation map.
Let $\psi_i$ 
be the first Chern class of the $i^{th}$ cotangent
line bundle $L_i$ on $\overline{M}_{0,n}(X, \beta)$. 
The fiber of $L_i$ over the 
moduli point $$[\mu: (C, p_1, \ldots, p_n)
\rarr X]$$
is the cotangent space of $C$ at $p_i$.
The Chern classes $\psi_i$
are elements in 
$H^2(\overline{M}_{0,n}(X, \beta),\mathbb{Q})$.
In [W2], invariants of $X$ are defined by integrals
over the moduli space of maps. The genus 0 {\em gravitational
descendents} are the invariants:
\begin{equation}
\label{wqed}
\lan \  \tau_{a_1}(\gamma_1) \cdots
\tau_{a_n}(\gamma_n)\ \ran_{\beta}^{X} =
\int_
{\overline{M}_{0,n}(X, \beta)} e_1^*(\gamma_1) \scup
\psi_1^{a_1} \scup \cdots \scup e_n^*(\gamma_n) \scup 
\psi_n^{a_n}
\end{equation}
where $\gamma_i\in H^*(X, \mathbb{Q})$ and the $a_i$ 
are nonnegative integers.
As usual, the invariants are defined to vanish unless the
dimension of the integrand is correct. When the $a_i$
are all 0, the gravitational descendents specialize to
the {\em Gromov-Witten invariants} of $X$. For simplicity,
$\tau_0(\gamma)$ will often be denoted by $\gamma$ in
(\ref{wqed}).

In this preliminary section, three topics are covered.
Basic properties of the descendent integrals are treated first.
Next, a fundamental solution of the
quantum differential equation obtained from the
flat connection in the Dubrovin formalism is derived.
Finally, an explicit form of this solution in case
$X=\proj^m$ is given.  The main sources in the mathematics
literature for this material are [D], [G1].

The formulas of fundamental class
and divisor for Gromov-Witten
invariants (see [RT], [KM]) take a slightly different form for
gravitational descendents. These formulas are given
in Section \ref{sdd} and are closely related
to Witten's equations in [W2]. 
Together with the topological recursion relations, 
these formulas are
sufficient to reconstruct the genus 0 descendents from
the Gromov-Witten invariants and to
derive a fundamental solution to the quantum differential
equation.

$X$ will be assumed to be convex. This hypothesis
leads to great simplification in the genus 0 case:
no virtual fundamental class considerations are needed at this point.
The results, however, are valid for general nonsingular
projective $X$.

\subsection{The string, dilaton,  and divisor equations}
\label{sdd}
Let the map $$\nu: \overline{M}_{0,n+1}(X, \beta)
\rarr \overline{M}_{0,n}(X, \beta)$$ be the natural
contraction morphism forgetting the last point. Contraction is
possible only when $\beta=0$, $n\geq 3$ or
$\beta\neq 0$, $n \geq 0$.
Three basic equations hold for descendent invariants:
the string, dilaton, and divisor equations.
They apply when contraction is possible and
the class assigned to the
last marking is of total codimension 0 or 1.

\begin{enumerate}
\item[I.] The string equation. Let $T_0 \in H^*(X, \Q)$ be the
unit:
$$
\lan \ \tau_{a_1}(\gamma_1) \cdots \tau_{a_n}(\gamma_n) 
\ T_0 \ \ran_\beta =$$
 $$\sum_{i=1}^n  \lan\ \tau_{a_1}(\gamma_1) \cdots 
\tau_{a_{i-1}}(\gamma_{i-1})\ 
\tau_{a_i-1}(\gamma_i)\ \tau_{a_{i+1}}(\gamma_{i+1})
\cdots \tau_{a_n}(\gamma_n)\ \ran_\beta.$$

\item[II.] The dilaton equation:
$$
\lan \ \tau_{a_1}(\gamma_1) \cdots \tau_{a_n}(\gamma_n) 
 \ \tau_{1}(T_0) \ \ran_\beta =
 (-2+n)\cdot  \lan \ \tau_{a_1}(\gamma_1) \cdots \tau_{a_n}(\gamma_n) 
\ \ran_\beta  $$

\item[III.] The divisor equation. Let $\gamma\in H^2(X, \mathbb{Q})$:
$$
\lan \ \tau_{a_1}(\gamma_1) \cdots \tau_{a_n}(\gamma_n)
 \ \gamma \ \ran_\beta = (\int_\beta \gamma) \ \cdot \ 
\lan \ \tau_{a_1}(\gamma_1) \cdots \tau_{a_n}(\gamma_n)
\ \ran_\beta $$
$$ + \sum_{i=1}^n  \lan \ \tau_{a_1}(\gamma_1) \cdots  
\tau_{a_{i-1}}(\gamma_{i-1}) \ \tau_{a_i-1}(\gamma_i\scup \gamma) 
 \ \tau_{a_{i+1}}(\gamma_{i+1})
\cdots \tau_{a_n}(\gamma_n) \ \ran_\beta$$
\end{enumerate}
In these formulas, any term with a negative exponent
on a cotangent line class is defined to be 0.

The proofs of these equations in the convex genus 0 case
rely on a comparison result for cotangent lines:
$\psi_i=\nu^* (\psi_i)+ [D_{i,n+1}]$. Here, $D_{i,n+1}$
denotes the boundary divisor in 
$\overline{M}_{0,n+1}(X, \beta)$ determined
by the curve and point partition (see [FP]): 
$$(\beta_1=0, \{i, n+1\}\ | \
\beta_2=\beta, \{1, \ldots,\hat{i}, \ldots, n\}).$$ Equations I-III follow
easily from this comparison result.
In the nonconvex or higher genus cases, the string and
divisor equation hold exactly in the form above.
The dilaton equation  
is true in genus $g$ with the factor $(-2+n)$
replaced by $(2g-2+n)$. The proofs in this greater generality
require properties of the virtual fundamental class.

In genus 0, the descendent integrals actually carry no
more information than the Gromov-Witten invariants (see [Du]):

\begin{pr}
\label{recc}
The genus 0 descendents of $X$ can be uniquely
reconstructed
from the genus 0
Gromov-Witten invariants.
\end{pr}

The proof is via the topological recursion relations.
For $n\geq 3$, consider the map
$$\nu: \overline{M}_{0,n}(X, \beta) \rarr \overline{M}_{0,3}$$
forgetting all data except the first 3 markings.
Again, a comparison result for cotangent lines is needed:
$\psi_1 - \nu^*(\psi_1)$ is seen to equal
a linear combination of boundary divisors of
$\overline{M}_{0,n}(X, \beta)$. Since $\psi_1$ is 0 in
$H^2(\overline{M}_{0,3}, \Q)$,
$\psi_1$ is a boundary class on 
$\overline{M}_{0,n}(X, \beta)$.
The divisors which occur in $\psi_1 - \nu^*(\psi_1)$
are those with point splitting $A\scup B$ where $1\in A$ and
$\{2,3\} \subset B$. No multiplicities occur (this may be seen, for example,
by intersections with curves). 

\begin{lm}
\label{trel}
Let $n\geq 3$. The following  boundary expression for $\psi$ holds:
$\psi_1 = \sum_{\Gamma} D$, 
where the sum is over all boundary divisors with
point splitting separating $1$ from $\{2,3\}$.
\end{lm}

Using Lemma \ref{trel}
together with the recursive structure of the boundary
(see [Ko], [BM], [FP]), a {\em topological
recursion relation} among genus 0 descendent integrals
is obtained. First, let $T_0, \ldots , T_m$ be a
basis of $H^*(X, \Q)$ (we assume here the cohomology is
all even dimensional to avoid signs). Let $g_{ef}= \int_X 
T_e \scup T_f$ be the intersection pairing, and
let $g^{ef}$ be the inverse matrix. The recursion relation is:
\begin{equation}
\label{trr}
\lan\ \tau_{a_1}(\gamma_1)\ \tau_{a_2}(\gamma_2)\
\tau_{a_3}(\gamma_3) \ \prod_{i\in S} \tau_{d_i}(\delta_i) \ \ran_\beta =
\end{equation}
$$ \sum  
\lan\ \tau_{a_1-1}(\gamma_1)\ \prod_{i\in S_1} \tau_{d_i}(\delta_i)\ T_e
\ \ran_{\beta_1} g^{ef}
\lan\ T_f \ \tau_{a_2}(\gamma_2)\ \tau_{a_3}(\gamma_3)\
 \prod_{i\in S_2} \tau_{d_i}(\delta_i)\ \ran_{\beta_2}.
$$
The sum is over all stable splittings $\beta_1+\beta_2= \beta$,
$S_1\scup S_2= S$, and over the diagonal splitting indices $e,f$:
the class $\sum g^{ef} T_e \otimes T_f$ is the Poincar\'e
dual of the diagonal $\deli \subset X \times X$.

The proof of the Reconstruction Theorem follows easily
from (\ref{trr}). The argument is via induction on
the number of cotangent line classes.  
A descendent with
no cotangent line classes is a Gromov-Witten invariant by
definition. All $\beta=0$ invariants are determined
by the classical cohomology of $X$ together with 
well-known formulas for cotangent line class
integrals on $\overline{M}_{0,n}$ [W2].
The topological recursion relations reduce descendents
with at least 3 markings to integrals with fewer 
cotangent line classes. Let $\lan I \ran_{\beta\neq 0}$
be a descendent integral with only 2 markings.
Let $H$ be an ample divisor on $X$.
Add an extra marking subject to the divisor $H$
condition: $\lan I \cdot H \ran_\beta$. The divisor equation then
relates 
$\lan I \ran _\beta$ and $\lan I\cdot H 
\ran _\beta$ modulo descendents with
fewer cotangent lines. Since $\lan I \cdot H \ran_\beta$ has 3 markings,
equation (\ref{trr}) equates $\lan I \cdot H \ran_\beta$  with an
expression involving descendents with fewer
cotangent lines. Similarly, if $\lan I \ran _{\beta\neq 0}$ is an
integral with only 1 marking, then $\lan I \cdot H  \cdot 
H \ran _\beta$
is considered. 
This completes the proof of Proposition \ref{recc}.

\subsection{The fundamental solution}
\label{fundsolly}
Gravitational descendents arise
in 
fundamental
solutions of the quantum differential equation of
$X$. 
Givental derives a solution naturally
via a related torus action and equivariant
considerations.
The solution is rederived here using the topological
recursion relations (as suggested to the author by
S. Katz). It is also possible to derive the solution
from the WDVV-equations for descendents. This solution
was found by Witten and Dijkgraaf; it also appears in [Du].

As before, let $T_0, \ldots, T_m$ be a homogeneous
basis of $V=H^*(X, \mathbb{Q})$ such that $T_0$ is the
ring identity and $T_m$ is its Poincar\'e dual.
The tangent space of $V$ at every point is canonically identified
with $V$. Let $\partial_0, \ldots, \partial_m$
be the corresponding tangent fields.
Let $\gamma=\sum t_iT_i$ be coordinates on $V$
defined by the basis. Let $F= \sum f^i \partial_i$ be
a vector field.
Let $\Phi$ be the quantum potential defined
by the Gromov-Witten invariants:
$$\Phi(t_0, \ldots,t_m)=\sum_{n\geq3} \sum_{\beta} \frac{1}{n!} \lan
\gamma^n\ran _\beta.$$
In the basic case $X=\mathbf{G}/\mathbf{P}$, 
the potential is a formal series in $\Q[[t_i]]$ 
(in general additional variables are needed).
The quantum product is defined by:
\begin{equation}
\label{bigpr}
\partial_i * \partial_j = \Phi_{ijr}g^{rs} \partial_s.
\end{equation}
Define a (formal) connection $\gradd$ on the tangent
bundle of $V$ by:
$$\graddi (F) =
\hb \frac{\partial F}{\partial t_i} - 
\partial_i * F=\sum \ (\hb \frac{\partial f^s}{\partial t_i}-
\Phi_{ijr} g^{rs}f^j) \ \partial_s.$$
The WDVV-equations imply that $\gradd$ is flat (see [D], [KM], [G1]).
Therefore, flat vector fields $F$ exist formally.
The equations for flat solutions $F$ are:
$\hb \partial F/\partial t_i = \partial_i * F$.
This is the quantum differential equation.

Following [G1], define a matrix of formal functions
in $\Q[[\hb^{-1}, t_i]]$:
$$\Psi_{ab}= g_{ab}+
\sum_{n\geq 0, \ \beta, \ (n,\beta)\neq (0,0)} \frac{1}{n!}
\lan T_a \cdot \frac{T_b}{\hb-\psi} \cdot \gamma^n\ran _\beta$$
where $0\leq a,b \leq m$.
The matrix may be written more explicitly as:
$$\Psi_{ab}= g_{ab}+
\sum \sum_{k\geq 0} \frac{\hb^{-k-1}}{n!} 
\lan T_a\cdot \tau_k(T_b) \cdot \gamma^n \ran _\beta$$
with the same summation conventions on $n$ and $\beta$.

\begin{pr}
$\Psi$ yields a fundamental solution of the quantum 
differential equation:
\begin{equation}
\label{pppy}
\graddi\  \sum \Psi_{ab}g^{as} \partial_s=0.
\end{equation}
The constant term of the solution $\Psi_{ab}g^{as}$
is the identity matrix.
\end{pr}

\bpf
By the definitions, it follows:
\begin{equation}
\label{fff}
\hb \frac{\partial \Psi_{ab}}{\partial t_i} g^{as}=
\sum_{n\geq 0,\ \beta}\ \sum_{k\geq 0} \frac{\hb^{-k}}{n!}
\lan T_a\cdot \tau_k(T_b) \cdot T_i \cdot \gamma^n \ran_\beta \ g^{as}.
\end{equation}
The coefficient of $\partial_s$
in the second term $\partial_i* \Psi_{ab}g^{aj}\partial_j$ 
in the covariant derivative is:
\begin{equation}
\label{sss}
\Phi_{ijr}g^{rs}\Psi_{ab}g^{aj}=
\Phi_{ibr}g^{rs}+
\end{equation}
$$\sum\sum_{k\geq 0} 
\frac{\hb^{-k-1}}{n_1!n_2!}
\lan T_a \cdot \tau_k(T_b) \cdot \gamma^{n_1}\ran_{\beta_1} \ g^{aj}
\lan T_j \cdot T_i \cdot T_r \cdot \gamma^{n_2} 
\ran_{\beta_2} \ g^{rs}$$
where the first sum is over stable splittings
$n_1+n_2=n$ , $\beta_1+\beta_2=\beta$,
(and, of course,
the repeated indices).

The $k=0$ terms of (\ref{fff}) sum to exactly the first
term $\Phi_{ibr}g^{rs}$ of (\ref{sss}).
The $k\geq 1$ terms of (\ref{fff}) may be replaced via
the topological recursion relations (\ref{trr}) relative to
the first 3 markings to obtain precisely the second term
in (\ref{sss}). \epf

The fundamental solution takes a simpler form when restricted
to the space $H^2(X, \Q)$ -- that is, after passing to the
small quantum cohomology.  
The {\em small quantum product} is obtained by setting
all variables $t_i$ to 0 in the formula (\ref{bigpr}) 
which do not correspond to 
cohomology basis elements in $H^2(X, \Q)$.
Let $T_1, \ldots, T_k$ span $H^2(X, \Q)$.
Let $T$ denote the vector of cohomology classes
$(T_1, \ldots, T_k)$.
Let $t$ denote the vector of variables $(t_1, \ldots, t_k)$.
For $\beta \in H_2(X, \Z)$, let $v_\beta$ denote
the vector of constants $(\int_\beta T_1, \ldots, \int_\beta T_k)$.
We assume the classes
$T_1, \ldots, T_k$  
pair {\em non-negatively} with all effective curve classes
in $X$ (this is possible for $X=\mathbf{G}/\mathbf{P}$, but
required only for simplicity). 
By the divisor
formula, the small product may be written as:
\begin{equation}
\label{smprod}
\partial_i * \partial_j = \sum_\beta 
e^{v_\beta \cdot t} \ \lan \ T_i\cdot  T_j \cdot
T_r \ \ran_\beta \  g^{rs}\ 
 \partial_s.
\end{equation}
The matrix $\Psi$ can be written 
after restriction to $H^2(X, \Q)$ as:
\begin{equation}
\label{smfund}
\Psi_{ab}= 
\sum_\beta
e ^{v_\beta \cdot t}\lan T_a \cdot \frac{e^{T\cdot t/\hb} T_b}
{\hb-\psi} \ran _\beta.
\end{equation}
Again, the divisor equation is used.
The formula (\ref{smfund}) is a sum of descendents
with 2 markings. The function of 
cohomology and cotangent classes at the second point is expanded to
define the invariant. A convention is made in (\ref{smfund})
regarding the $\beta=0$ case (as 2 point degree 0
invariants are not defined):
$$\lan T_a \cdot \frac{e^{T\cdot t/\hb} T_b}
{\hb-\psi} \ran _0 = \lan T_a \cdot  e^{T\cdot t/\hb} T_b 
\cdot 1\ran_0.$$
The series (\ref{smfund}) is viewed as a formal
power series in the variables $\hb^{-1}, t_i,$ and $e^{t_i}$. More
precisely, the series is an element of 
$\Q[\hb^{-1},t_i][[e^{t_i}]]$. Modulo the variables
$t_i$ and $e^{t_i}$, 
$\Psi_{ab} g^{as}$ is the constant identity matrix.

The small quantum differential equation is:
$$1\leq i \leq k, \ \ \ 
\hb \frac{\partial F}{\partial t_i} = \partial_i * F$$
where $F$ is a vector field function of only $t$, and the
product is the small quantum product.
For $1\leq i \leq k$, the small analogue of (\ref{pppy})
holds for (\ref{smfund}):
$$\hb \frac{\partial \Psi_{ab}}{\partial t_i} g^{as} \partial_s =
\partial_i * \Psi_{ab} g^{as} \partial_s.$$
In fact, the restricted matrix $\Psi_{ab}g^{as}$ provides
a fundamental solution to this small quantum differential
equation. 
Only the small quantum objects
will be considered in this paper.

\subsection{Projective space}
We now let $X=\proj^m$. Let $H$ denote the hyperplane
class in $H^2(\proj^m, \Q)$. Let $T_i= H^i$ be the
cohomology basis. 
Let $t=t_1$.
The small quantum ring structure is:
$$QH^*_{s} \proj^m = \Q[\partial_1, e^{t}]/ (\partial_1^{m+1}- e^{t})$$
(see [G1], [FP]).
Let $\sum_{0}^{m} f^i \partial _i$ be a vector field 
where $f^i=f^i(t)$.
The small quantum differential equation
is then the following system:
$$i>0, \ \ \hb \frac{\partial f^i}{\partial t} = f^{i-1}$$
$$  \hb \frac{\partial f^0}{\partial t} = e^{t} f^m.$$
The function $f^m$ determines a vector solution
if and only if it is annihilated by the operator
$\mathcal{D}=(\hb d/dt)^{m+1} - e^{t}$.
A (formal) fundamental solution to the equation $\mathcal{D}f=0$
is given by the following expression:
\begin{equation}
\label{solly}
S=S_{\proj^m} = \sum_{d\geq 0} \frac{e^{(H/\hb+d)t}}{ \prod_{r=1}^{d}
(H+r\hb)^{m+1}} \ \ mod \ H^{m+1}.
\end{equation}
$S$ is expanded in powers of $H$ (subject to $H^{m+1}=0$) as:
$$S = \sum_{b=0}^{m} S_b H^{m-b} $$
where $S_{b}$ is a formal series in $\Q[\hb^{-1},t][[ e^{t}]]$.
It is easily checked by formula (\ref{solly}) that
$\mathcal{D}$ annihilates $S_{b}$. Define the matrix $M$ of
functions by
$$M^s_{b}= (\hb \frac{d}{dt})^{m-s} S_b.$$
Modulo $t$ and $e^{t}$, the only contribution
to $M^s_b$ occurs in
the
$d=0$ summand in (\ref{solly}); it is the identity matrix.
 $M^s_{b} \partial _s $ defines a fundamental
solution to the small quantum differential equation.
By uniqueness,
\begin{equation}
\label{eqqq}
\Psi_{ab} g^{as} = M^s_b.
\end{equation}
The required uniqueness statement here depends on the
equality modulo $t, e^{t}$ and the fact that the
solutions lie in $\Q[\hb^{-1},t][[ e^{t}]]$.

Consider the hyperplane embedding
$X=\proj^m \subset \proj^{m+1}$.
The solution $S_{\proj^m}(t,\hb=1)$ 
agrees with the 
hypergeometric series $S^*_{\proj^m}$ defined in equation (\ref{hgs}) 
for the hyperplane $X \subset \proj^{m+1}$. This is the
first example of Givental's correspondence (i) of Section 0.

Equations (\ref{solly}) and (\ref{eqqq}) together with the
solution (\ref{smfund}) compute  
all 2 point invariants of $\proj^m$ with a cotangent
line class on 1 point. For example, 
tracing through the equations yields:
$$\lan\  \tau_{dm+d-2}(T_m)\ \ran_d =
 \int_{\overline{M}_{0,1}(\proj^m,d)}
e_1^*(H^m) \scup \psi^{dm+d-2} = \frac{1}{(d!)^{m+1}}.$$
The solution to the small quantum differential equation
provides an elegant organization of these 2 point
descendents.

For a general space $X$, define
$S_b=\Psi_{0b}$.
Let $\mathcal{D}(\hb, \hb \partial/ \partial t_i, e^{t_i})$
be a differential operator which is a
polynomial in the operators
$\hb, \hb \partial/ \partial t_i$, and $e^{t_i}$.

\begin{lm}
\label{van}
If $\mathcal{D}  S_b =0$ for all $0 \leq b \leq m$, then
the equation $\mathcal{D}(0, \partial_i, e^{t_i})=0$ holds
in $QH_{s}^*(X)$.
\end{lm}

\bpf
Let $M= \Psi_{ab} g^{as}$ be the solution matrix.
Let $p=\mathcal{D}(0, \partial_i, e^{t_i})$ in
$QH_{s}^*(X)$. Let $P$ be the matrix of quantum
multiplication by $p$ in the basis $\partial_0, \ldots, \partial_m$.
By the quantum differential equation,
\begin{equation}
\label{ddd}
\mathcal{D} M  =  
P \cdot  M  + \sum_{k\geq 1}^K \hb^k C_k \cdot M
\end{equation}
where $C_k$ is a matrix with coefficients in $\Q[[e^{t_i}]]$ and
$K<\infty$.
$M$ is certainly invertible in the coefficient
ring  $\Q[[\hb^{-1}, t_i, e^{t_i}]]$. 
The vector
$(S_0, \ldots, S_m)$ is the $m^{th}$ row of the
solution matrix $M$.  
Hence, after multiplying (\ref{ddd}) by $M^{-1}$ from the
right and using the $\hb$-grading, it follows that the $m^{th}$ row of $P$
vanishes.   Poincar\'e duality on $X$ and
the definition of the quantum product 
then imply $p=0$. \epf

\section{\bf The $*_X$-product induced by a hypersurface}
\subsection{The product construction}
Let $Y=\proj^m$ (or more generally a homogeneous
variety $Y=\mathbf{G}/\mathbf{P}$). 
Let $L$ be an ample line bundle on $Y$. Let $s\in H^0(Y,L)$
be a section with nonsingular zero locus $X$. 
Assume the
dimension of $X$ is at least 3. The isomorphisms
\begin{equation}
\label{idds}
H_2(Y,\Z)\eqq H_2(X, \Z), \ \ H^2(Y,\Q) \eqq H^2(X,\Q)
\end{equation}
hold by the Lefschetz theorem.

Gromov-Witten invariants are defined on $X$ via the
virtual fundamental class since the moduli space
of maps is in general ill-behaved. In our special
situation, the virtual class on $\overline{M}_{0,0}(X,\beta)$
has a simple interpretation.
Let $\pi:U \rarr \overline{M}_{0,0}(Y,\beta)$ be the
universal curve. Let $\mu: U \rarr Y$ be the universal
map. It is not hard to check the following facts.
The sheaf $E_\beta=\pi_* \mu^*(L)$ is a vector bundle
of rank $\int_\beta c_1(L) + 1$ equipped with
a canonical section $s_E$ obtained from $s$. The moduli
space $$\overline{M}_{0,0}(X,\beta) \subset \overline{M}_{0,0}
(Y, \beta)$$ is the stack
theoretic zero locus of $s_E$. The virtual fundamental
class of $\overline{M}_{0,0}(X,\beta)$ is simply
the refined top Chern class of $E_\beta$ with respect to
the section $s_E$: it is a Chow class on 
$\overline{M}_{0,0}(X,\beta)$ of expected dimension.
Let $\nu: \overline{M}_{0,n}(Y, \beta) \rarr
\overline{M}_{0,0}(Y,\beta)$ be the forgetful morphism.
The virtual class of $\overline{M}_{0,n}(X,\beta)$
is the refined top Chern class of the canonical pull-backs of
$E_\beta$ and $s_E$ via $\nu$. For notational
convenience, the $\nu$ pull-back of $E_\beta$ will 
also be denoted $E_\beta$. The algebraic 
theory of the virtual fundamental class is developed
in [LT], [B], [BF]. The above construction appears in
[Ko]. 

Let $T_0, \ldots, T_m$ be a basis of $H^*(Y,\Q)$
satisfying the conventions of Section \ref{fundsolly}. 
As before, we assume that 
$T_1, \ldots, T_k$ is a basis of $H^2(Y, \Q)$ 
pairing {\em non-negatively} with all effective curve classes
in $Y$. 
Let $t_1, \ldots, t_k$
be the corresponding divisor variables, let $q_i= e^{t_i}$, and let
$q$ denote $(q_1, \ldots, q_k)$.
Let $QH^*_{s}(X) = H^*(X, \Q) \otimes \Q[[q]]$
with the canonical free $\Q[[q ]]$-module structure.
The small quantum product on $QH^*_{s}(X)$ is 
$\Q[[q ]]$-linear and 
defined via the
3 point genus 0 invariants as in (\ref{smprod}). 
Associativity is established
via properties of the virtual fundamental class [LT], [B].

A new $\Q[[q]]$-linear small quantum product $*_X$ is now defined 
on the free module $H^*(Y,\Q)\otimes \Q[[q]]$. Let
$i: X \rarr Y$ denote the inclusion.
Let 
\begin{equation}
\label{mapback}
i^* : H^*(Y, \Q) \otimes \Q[[q]] 
\rarr H^*(X, \Q) \otimes \Q[[q]]
\end{equation}
denote the canonical $\Q[[q]]$-linear pull-back.
The new product $*_X$ on the right side of (\ref{mapback})
will satisfy a homomorphism property:
for  $a, b \in H^*(Y, \Q)$,
\begin{equation}
\label{comppy}
 i^*( a *_X b) = (i^* a) * (i^* b)
\end{equation}
where the right side of (\ref{comppy}) is the product
in $QH^*_{s}(X)$.

The product $*_X$ is defined by new integrals over
the moduli space of maps to $Y$.
Consider again the vector bundle $E_\beta$ on
$\overline{M}_{0,n}(Y,\beta)$. For each 
marking $j$, there is a canonical bundle sequence
obtained by evaluation:
\begin{equation}
\label{bsq}
0 \rarr E'_{\beta,j} 
\rarr  E_{\beta} \rarr e_{j}^*(L) \rarr 0
\end{equation}
The bundle $E'_{\beta,j}$ is of rank $\int_\beta c_1(L)$.
If $\beta=0$, then $E'_{\beta,j}=0$.
Define new integrals: for 
cohomology classes $\gamma_1, \ldots, \gamma_{n} \in H^*(Y, \Q)$,
\begin{equation}
\label{inty}
\lan \gamma_1 \cdots \gamma_{n-1} \cdot \tilde{\gamma}_{n} 
\ran_{\beta}^{Y/X}= \int_{\overline{M}_{0,n}
(Y, \beta)} \prod_{i=1}^{n} e_i^*(\gamma_i) 
\scup c_{\text{top}}(E'_{\beta,n}).
\end{equation}
The tilde in the argument of (\ref{inty}) denotes the
marking with respect to which the construction (\ref{bsq})
is undertaken.
Define the new product $a*_X b$ for $a,b \in H^*(Y, \Q)$ by:
$$a *_X b = \sum_{\beta} 
q^{v_\beta}
\lan a \cdot b \cdot \tilde{T}_e\ran_{\beta}^{Y/X}
    g^{ef} T_f$$
where $q^{v_\beta}= q_1^{\int_\beta T_1} \cdots q_k^{\int_\beta T_k}$.

\begin{pr}
The product $*_X$ defines a commutative, associative, unital
ring structure on $H^*(Y, \Q) \otimes \Q[[q]]$ (with unit $T_0$).
\end{pr}

\bpf
The product is commutative by the symmetry of the
integrals (\ref{inty}) in the first two factors. 
The unital property of $T_0$ follows
for exactly the same reasons as in the usual quantum
product: only degree 0 terms contribute to a product
with $T_0$ and $\lan T_0 \cdot \gamma \cdot  \tilde{T}_e \ran_0^{Y/X} =
\int_Y \gamma \scup T_e$.
 
In the
usual quantum product, associativity is a consequence of
the basic boundary linear equivalence on $\overline{M}_{0,4}$
pulled back to $\overline{M}_{0,4} (Y, \beta)$.
A slight twist is needed here.
Let $a,b,c, \gamma \in H^*(Y, \Q)$.
Let $D(1,2\ |\ 3, 4)$ and $D(1,4\ | \ 2,3)$
be the divisors 
on $\overline{M}_{0,4}(Y, \beta)$
determined by the associated point splittings
(see [FP]). 
Let $$\omega =e_1^*(a) \scup e_2^*(b) \scup
e_3^*(c) \scup e_4^*(\gamma) \scup c_{\text{top}}(E'_{\beta,4})$$
There is an equality:
$$\int_{D(1,2\ |\ 3,4 )} \omega =
\int_{D(1,4\ |\ 2,3)} \omega.$$
The recursive structure of the boundary and the
simple behavior of the restriction of $E'_{\beta,4}$ yields:
$$\int_{D(1,2\ |\ 3,4 )} \omega =
\sum_{\beta_1+\beta_2=\beta} \lan a\cdot b \cdot \tilde{T}_e \ran
_{\beta_1}^{Y/X}
g^{ef} 
\lan T_f \cdot c \cdot\tilde{\gamma} \ran_{\beta_2}^{Y/X},$$
$$\int_{D(1,4\ |\ 2,3 )} \omega =
\sum_{\beta_1+\beta_2=\beta} \lan a\cdot \tilde{\gamma}\cdot T_e \ran
_{\beta_1}^{Y/X}
g^{ef} 
\lan \tilde{T}_f \cdot b \cdot c \ran_{\beta_2}^{Y/X}.$$
Associativity now follows easily. \epf

The following Proposition is due to T. Graber.

\begin{pr}
\label{compat}
The pull-back $i^*$ is a ring homomorphism from
$$QH_{s}^*(Y/X)=(H^*(Y, \Q) \otimes \Q[[q]], *_X)$$ to $QH_{s}^*(X)$.
\end{pr}

\bpf
Let $a,b\in H^*(Y, \Q)$. Let $a*_X b= \sum_{\beta} q^{v_\beta} c_\beta$.
Let the product of the pull-backs be:
$i^*(a) * i^*(b)= \sum_{\beta} q^{v_\beta} c'_\beta$.
The equality $i^*(c_\beta)= c'_\beta$ must be proven.

Consider the following fiber square:
\begin{equation*}
\begin{CD}
 Z  @>>>  \overline{M}_{0,3}(Y, \beta) \\
@V{e_3}VV   @V{e_3}VV \\
 X @>{i}>>  Y
\end{CD}
\end{equation*}
where $Z$ is the zero locus of $e_3^*(s) \in H^0(e_3^*(L))$.
Let $\gamma= e_1^*(a) \scup e_2^*(b)$ in $H^*(\overline{M}_{0,3}
(Y,\beta), \Q)$. Then,
$$c_\beta= e_{3*}(\gamma \scup
c_{\text{top}}(E'_{\beta,3})).$$
By properties of the Gysin map [F],
\begin{equation}
\label{ffff}
i^*(c_\beta)= e_{3*} i^! (\gamma \scup
c_{\text{top}}(E'_{\beta,3})).
\end{equation}
Recall the embedding
\begin{equation}
\label{qwq}
\overline{M}_{0,3}(X,\beta) \subset \overline{M}_{0,3}(Y,\beta)
\end{equation}
as the zero section of $E_{\beta}$.
By the realization of the virtual fundamental class of
maps to $X$,
\begin{equation}
\label{ssss}
c'_\beta= e_{3*} 0_{E_\beta}^! (\gamma).
\end{equation}
The inclusion (\ref{qwq}) factors
through $Z$. There is an equality of
classes on $Z$:
$$ 0_{E_\beta}^!(\gamma) = 0_{e_3^*(L)}^!(\gamma \scup c_{\text{top}}
(E'_{\beta,3}))
=
i^! (\gamma \scup c_{\text{top}}(E'_{\beta,3})).$$
The first equality follow from  sequence (\ref{bsq}); the
second is by definition.
Equation (\ref{ffff}) and (\ref{ssss}) yield the equality
$i^*(c_\beta)= c'_\beta$,
which concludes the proof. \epf

The small quantum differential equation for the product $*_X$
may also be considered:
\begin{equation}
\label{smcon}
\hb \frac{\partial}{\partial t_i} F = \partial_i *_X F.
\end{equation}
Again the coordinates $t=(t_1, \ldots, t_k)$ 
correspond to the cohomology
basis of $H^2(Y, \Z)$ and 
 $F$ is a vector field function of $t$.
The fundamental solution takes a form similar to (\ref{smfund}).
Let
\begin{equation}
\label{funcon}
\Psi_{ab} = \sum_{\beta} e^{v_\beta \cdot t} 
\lan \tilde{T}_a \cdot 
\frac{e^{T\cdot t/\hb} T_b}{\hb- \psi} \ran^{Y/X}_\beta.
\end{equation}
Again, a convention is required for $\beta=0$:
$$ \lan \tilde{T}_a \cdot 
\frac{e^{T\cdot t/\hb} T_b}{\hb- \psi} \ran^{Y/X}_0 =
\lan \tilde{T}_a \cdot 
{e^{T\cdot t/\hb} T_b} \cdot 1 \ran^{Y/X}_0 =
\lan T_a \cdot 
{e^{T\cdot t/\hb} T_b} \cdot 1 \ran^{Y}_0.$$
The fundamental solution of (\ref{smcon}) is 
$\Psi_{ab} g^{as} \partial _s$. This may be proven as
a specialization of the fundamental solution
$$\Psi_{ab}= g_{ab}+
\sum_{n\geq 0, \ \beta, \ (n,\beta)\neq (0,0)} \frac{1}{n!}
\lan \tilde{T}_a \cdot \frac{T_b}{\hb-\psi} \cdot \gamma^n\ran _\beta$$
of the differential equation obtained from large product $*_X$ 
$$\partial_i *_X \partial_j =
\sum_{n\geq 0} \sum_\beta 
\frac{1}{n!}\lan \ T_i \cdot T_j \cdot \tilde{T}_r \cdot \gamma^n \ 
\ran_\beta^{Y/X}g^{rs}\partial_s , \ \ \
\gamma= \sum_{i=0}^m t_i T_i$$
via 
analogous topological recursion
relations for the tilde integrals  (following the
proofs in Section 1). Alternatively, (\ref{funcon}) may be
proven to be a fundamental solution directly from the
definition (\ref{smcon}) together with 
the tilde topological recursion relations.

\subsection{Projective space}
\label{vbvbvb}
Let $Y=\proj^m$. Let $X$ be a hypersurface of degree $l$.
Consider the fundamental solution of (\ref{smcon}) 
given by the matrix $\Psi_{ab}g^{as}$. Let $t=t_1$.
Let $S_{b}= \Psi_{ab}g^{a, m-1}$ be the $(m-1)^{st}$ row 
of the solution matrix. A Lemma analogous to 
Lemma \ref{van} holds.
Let $\mathcal{D}(\hb, \hb \partial/ \partial t, e^{t})$
be a differential operator which is a
polynomial in the operators
$\hb, \hb \partial/ \partial t$, and $e^{t}$.
\begin{lm}
\label{wwww}
If $\mathcal{D}  S_b =0$ for all $0 \leq b \leq m$, then
the equation $\mathcal{D}(0, \partial_1, e^{t})=0$ holds
in $QH_{s}^*(X)$.
\end{lm}

\bpf
Let $p=\mathcal{D}(0, \partial_1, e^{t})$ in
$QH_{s}^*(Y/X)$. Let $P$ be the matrix of quantum 
$*_X$-multiplication by $p$ in the basis $\partial_0, \ldots, \partial_m$.
By the argument used in the proof of Lemma \ref{van},
it follows that the $(m-1)^{st}$ row of $P$
vanishes.   
The definition of the product $*_X$ together with
classical intersection theory on $Y$
then implies $p= f(e^t) \partial_m$. Hence, $i^*(p)=0$.
The Lemma then follows from Proposition \ref{compat}.\epf

The functions $S_b$ will be calculated 
for hypersurfaces of degree $l \leq m+1$ via
torus localization. For notational simplicity, let
$H^b=T_b$.
We first organize the functions $S_b$ in a more convenient form: 
$$\tilde{S}(t, \hb)= \sum_{b=0}^m 
S_b H^{m-b} \in H^*(Y, \Q)[\hb^{-1},t]
[[e^t]].$$
The definitions imply: 
\begin{eqnarray}
\tilde{S}(t,\hb) & = & \sum_{b=0}^m \sum_{d\geq 0} e^{d t} 
\lan \tilde{H} \cdot 
\frac{e^{H t/\hb} T_b}{\hb- \psi} \ran^{Y/X}_d  H^{m-b}\\
\label{fredd}
& = & \sum_{d\geq 0} e^{(H/\hb+d) t}   e_{2*}( 
\frac{1}{l} \cdot
\frac {c_{\text{top}}(E_d)}{\hb-\psi_2}). 
\label{frod}
\end{eqnarray}
In degree 0, the convention
\begin{equation}
\label{deggo}
e_{2*}( \frac{c_{\text{top}}(E_d)}{\hb-\psi_2}) =l \cdot  H \in H^*(Y, \Q)
\end{equation}
is taken. 

\begin{lm}
$\tilde{S}(t,\hb)$ is divisible by $H$.
\end{lm}
\bpf
The degree 0 contribution is clearly divisible by $H$
by (\ref{deggo}).
It suffices to show  the integrals 
$$\langle\tilde{H} \cdot \frac{T_m}{\hb-\psi}\rangle^{Y/X}_d= 
\int_{[\overline{M}_{0,2}(X,d)]^{\text{vir}}} 
\frac{1}{l}\cdot
\frac{e_2^*(T_m)}
{\hb-\psi_2}$$
vanish for $d>0$. 
In fact, these are zero for a trivial reason: the point
class $T_m$ does not intersect $X$ in $Y$. \epf

The correlator $S(t, \hb)= l\cdot \tilde{S}(t,\hb)$ is more
convenient for the calculations to come.
The dependence on $l$ will be made
explicit as a third argument $S(t,\hb,l)$.
Givental's correlator $S_X$ for
a hypersurface $X \subset 
\proj^m$ is defined by 
\begin{equation}
\label{deffff}
 S_X(t,\hb) =\frac{1}{lH}\  S(t, \hb,l).
\end{equation}
The correlators $S_X$ considered in the correspondences (i)-(iii) of
Section 0 are evaluated at $\hb=1$.
It is $S(t, \hb,l)$ which 
is explicitly calculated in Sections 3 and 4 via fixed
point localization for the natural torus action on the moduli
space of maps.

\section{\bf Hypersurfaces: Cases (i) and (ii)}
\subsection{The $\mathbf{T}$-action}
\label{aaction}
We first set notation for the required torus action.
Let $\lambda$ denote the set $\{\lambda_0, \ldots, \lambda_m\}$.
Let $\R=\Q[\lambda]$ 
be the standard presentation of the equivariant 
cohomology ring of the torus $\mathbf{T}=(\com^*)^{m+1}$ :
$\lambda_i$ is the equivariant first Chern class of the 
dual of the standard
representation of the $i^{th}$
factor $\com^*$.
Let $\mathbf{T}$ act on $W=\com^{m+1}$ 
via the
standard diagonal representation.
A $\mathbf{T}$-action on $\proj(W)=\proj^{m}$ is canonically obtained.
The $\mathbf{T}$-action lifts canonically to $\oh_{\proj^m}(1)$
(and therefore to $\oh_{\proj^m}(k)$ for every $k$).
Let $H \in H^*_{\mathbf{T}}(\proj^m)$ be the
equivariant first Chern class of $\oh_{\proj^m}(1)$. All
equivariant cohomology groups will be taken with $\Q$-coefficients.
The standard presentation of $
H^*_{\mathbf{T}}(\proj^m)$ is:
\begin{equation}
\label{prez}
H^*_{\mathbf{T}}(\proj^m) \eqq  \Q[H, \lambda]/
(\Pi_{i=0}^m
(H-\lambda_i)).
\end{equation}
The fixed points $\{p_0, \ldots, p_m\}$ 
of the $\mathbf{T}$-action on $\proj^m$
correspond to the basis vectors in $\com^{m+1}$.
Let $\phi_i \in H^m_{\mathbf{T}}(\proj^m)$ denote (the dual of) the
equivariant fundamental class of the point $p_i$.
There is a naturally graded equivariant push-forward:
 $$\int_{\proj^m}: H^*_{\mathbf{T}}(\proj^m) \rarr \R.$$
For $x,y\in H^*_{\mathbf{T}}(\proj^m)$,
an inner product $\langle x, y \rangle = \int_{\proj^m} x\scup y$
is defined. The elements $\phi_i$ satisfy the following
basis property:
\begin{equation}
\label{vbvb}
x=y \ \ \ \Longleftrightarrow \ \ \ \forall i,\   
\langle \phi_i, x \rangle = \langle \phi_i, y\rangle.
\end{equation}
Property (\ref{vbvb}) follows easily from the
presentation (\ref{prez}).

There is a canonically induced $\mathbf{T}$-action on the 
stack
$\overline{M}_d= \overline{M}_{0,2}(\proj^m,d)$.
The torus acts on a stable
map to $\proj^m$ by translating the image.
The $\mathbf{T}$-fixed locus has been
determined in [Ko], and is (necessarily) a nonsingular substack. 
Each domain component of a $\mathbf{T}$-fixed
map must have as image in $\proj^m$ a $\mathbf{T}$-orbit of dimension
0 or 1. The 0-dimensional orbits are the
$m+1$ fixed points, and the 1-dimensional orbits are
the $\binom{m+1}{2}$ lines connecting these points. Moreover,
for a map to be $\mathbf{T}$-fixed, all nodes, marking, and ramifications
point must have fixed images in $\proj^m$.
As a result, the components of the 
$\mathbf{T}$-fixed locus of $\overline{M}_d$
are in bijective correspondence with graph types describing the
configurations of the domain components and
the markings of the map. 

More precisely, the graphs arising
in this correspondence are triples $(\Gamma, \mu, \delta)$
where $\Gamma$ is $2$-pointed tree, 
$$\mu: \text{Vert}(\Gamma)\rarr
\{p_0, \ldots, p_m\}, \ \ \text{and} \ \ \delta: \text{Edge}(\Gamma) \rarr 
\mathbb{Z}^{> 0}.$$
Let $\zeta=[\mu:(C,x_1,x_2) \rarr \proj^m] \in \overline{M}_d^{\mathbf{T}}$.
The graph $(\Gamma, \mu, \delta)$ associated to the fixed
component containing $\zeta$ is constructed as follows. 
The vertices $v\in\text{Vert}(\Gamma)$ correspond to 
 {\em connected components} $D(v)$
of $\mu^{-1}(\{p_0, \ldots, p_m\})$. 
$D(v)$ may be 0 or 1 dimensional.
The map
$\mu$ on $\text{Vert}(\Gamma)$ is determined by:  
$\mu(v)= \mu(D(v))$. 
An edge $e$ connecting vertices $v,v'$ corresponds to a
component $D(e)=\proj^1\subset C$  incident to $D(v), D(v')$ and 
lying over the $\mathbf{T}$-invariant line $L$
connecting the fixed points $p_\mu(v), p_\mu(v')$.
The degree of the map from $D(e)$ to $L$ is $\delta(e)$
(this
map is uniquely determined (up to isomorphism) by $\delta(e)$
since it is unramified over
$L\setminus \{p_\mu(v), p_\mu(v')\}$).
The markings on $C$ determine vertex markings on $\Gamma$. 
We will let the symbol $\Gamma$ denote the entire decorated graph structure.
Let 
$$\widetilde    {M}_\Gamma =
\Pi_{v\in \text{Vert}(\Gamma)}\  \overline{M}_{0,\text {val}(v)},$$
where $\overline{M}_{0,1}$ and $\overline{M}_{0,2}$ are taken to be points.
The valence of a vertex includes the markings as well as the incident edges.
The fixed component associated to
a graph $\Gamma$ is simply the stack quotient: 
$\overline{M}_\Gamma =
\widetilde{M}_\Gamma /
\mathbf{G}$,
where $\mathbf{G}$ is an associated automorphism group (see [GP]).
The order of $\mathbf{G}$ is equal to 
$$\Pi_{e\in \text{Edge}(\Gamma)}\delta(e) \cdot |\text{Aut}(\Gamma)|,$$
where $\text{Aut}(\Gamma)$ is the automorphism group of the
decorated graph.
An excellent description of these graphs may be found in [Ko].

Let $X\subset \proj^m$ be a hypersurface of degree $l\leq m+1$.
We follow the notation of Section 2 (with the
$l$ dependence explicit in the correlators).
As all structures in the definition of $S(t, \hb,l)$ are
canonically $\mathbf{T}$-equivariant, an element $S_{\mathbf {T}}(t,\hb,l)
\in H^*_{\mathbf{T}}(\proj^m, \Q)[[\hb^{-1}, t,e^t]]$
is defined by:
\begin{equation}
\label{lindad} 
S_{\mathbf{T}}(t,\hb,l)= \sum_{d\geq 0} e^{(H/\hb+d) t}   e_{2*}( 
\frac {c_{\text{top}}(E_d)}{\hb-\psi_2}).
\end{equation}
In the following definitions and
computations, all the geometric 
structures (schemes, sheaves, maps, push-forwards,
cohomology groups) will be given their canonical 
$\mathbf{T}$-equivariant interpretations. 

\subsection{Linear recursions in cases (i) and (ii)}
The calculation of $S_{\mathbf{T}}(t,\hb,l)$ relies on linear
recursions obtained from localization.
The recursions involve a related set of equivariant 
correlators $Z_i \in 
\R[[\hb^{-1},q]]$
defined for $0\leq i \leq m$ by:
\begin{equation}
\label{zdef}
Z_i(q, \hb,l) = 1+ \sum_{d> 0}  q^d \int_{\overline{M}_d} 
\frac{E'_{d,2}}{\hb-\psi_2}e_2^*(\phi_i).
\end{equation}
The integral on the right is the equivariant push-forward to 
a point. $E'_{d,2}$ in the integrand denotes the top Chern class 
of the bundle $E'_{d,2}$. The latter convention will be kept throughout
this section:
bundles appearing in integrands will always denote their
top Chern classes.
By the definitions of $E'_{d,2}$, $S_{\mathbf{T}}(t, \hb,l)$, and
the exact sequence (\ref{bsq}), we see:
\begin{equation}
\label{prrod}
\langle \phi_i, S_{\mathbf{T}}(t,\hb,l) \rangle = 
e^{\lambda_i t/\hb} l\lambda_i
Z_i(e^t, \hb,l),
\end{equation}
where the pairing $\langle, \rangle$ is
taken to be linear in the auxiliary parameters $\hb^{-1},t, e^t$.
The degree 0 term involves a matching of conventions.
By equation (\ref{prrod})
and property (\ref{vbvb}), 
the correlators $Z_i$ determine $S_{\mathbf{T}}$.

The dimension of $\overline{M}_{d}$ is $(m+1)d+m-1$,
the rank of $E'_{d,2}$ is $ld$, and the codimension of
the class $e_2^*(\phi)$ is $m$. Therefore, for $l\leq m$,
initial terms
in the $(1/\hb)$-expansion of (\ref{zdef}) 
with $\psi_2$-degree less than $(m+1-l)d-1$
vanish by dimension
considerations. Hence, we find:
\begin{equation}
\label{hobo}
Z_i(q, \hb,l\leq m) = 1+ \sum_{d > 0}  \Big(\frac{q}{\hb^{m+1-l}}\Big)^d 
\int_{\overline{M}_d} 
\frac{\psi_2^{(m+1-l)d-1}}{1-\psi_2/\hb} E'_{d,2} e_2^*(\phi_i).
\end{equation}
We rescale the $q$-dependence in these cases by the following definition:
$z_i(Q,\hb,l\leq m)=Z_i(Q\hb^{m+1-l}, \hb,l)$. Such vanishing
does not apply in the Calabi-Yau case.

Equivariant integrals over $\overline{M}_d$ may be computed 
via $\mathbf{T}$-equivariant localization.
Let $I$ be an equivariant class in
$H^*_{\mathbf{T}}(\overline{M}_d)$.
Let $G_d$ denote the 
set of components of the $\mathbf{T}$-fixed stack 
$\overline{M}_d^{\mathbf{T}}$.
The elements of $G_d$ are labelled by decorated graphs $\Gamma$. 
The equivariant integral of $I$ over $\overline{M}_d$ equals a sum
of contributions over the set $G_d$:
\begin{equation}
\label{bbott}
\int_{\overline{M}_d} I = \sum _{\Gamma\in G_d}
\int_{\overline{M}_\Gamma} \frac{I}{N_\Gamma},
\end{equation}
where $N_\Gamma$ is the equivariant normal bundle to
$\overline{M}_\Gamma \subset \overline{M}_d$.
Equation (\ref{bbott}) is the Bott residue formula
in equivariant cohomology. Explicit formulas
for the equivariant Euler class of $N_\Gamma$ in terms
of tautological classes in $H^*_{\mathbf{T}} (\overline{M}_\Gamma)$
have been obtained by Kontsevich in [Ko] (see also [GP], [CK]).
Givental requires the full formulas four times.
In each case, the formulas are
used to verify an algebraic equality -- 
the method is 
straightforward algebraic manipulation.
The computations in the Calabi-Yau case will be covered in
some detail. The others will be omitted.

Let $d>0$.
Fix an index $0\leq i \leq m$, and consider the
integrals in (\ref{zdef}) and (\ref{hobo}). In order to
analyze these integrals, it is necessary to study the
graphs associated to the fixed loci.
We partition the set $G_d$ into 3 disjoint subsets:
$$G_d= G_d^{i*} \cup G_d^{i0} \cup G_d^{i1}.$$
The set $G_d^{i*}$ consists of the fixed loci for which the
$2^{nd}$ marked point on the domain curve is {\em not} mapped to 
$p_i\in \proj^m$. 
The set $G_d^{i0}$ consists of
loci for which an {\em irreducible component} of the domain curve containing
the $2^{nd}$ marked point is collapsed to 
$p_i$. Finally, $G_d^{i1}$  consists of
loci for which the $2^{nd}$ marking is mapped to 
$p_i$ without lying on a collapsed
component. 
Let $G^i_d= G^{i0}_d \cup G^{i1}_d$.
Let $G^{i0}$ and $G^{i1}$ denote the
unions $\bigcup_{d> 0} G^{i0}_d$ and $\bigcup_{d>0} G_d^{i1}$
respectively. The three graph types have the following
basic properties:  
\vspace{+12pt}

\noindent{\bf{Type} $G_d^{i*}$.} Let $\Gamma \in G_d^{i*}$. 
As $e_2^*(\phi_i)$
vanishes when restricted to $\overline{M}_\Gamma$, the contribution of
$\Gamma$
to the integrals in (\ref{zdef}) and (\ref{hobo}) via
(\ref{bbott}) is 0.
\vspace{+12pt}

\noindent{\bf{Type} $G_d^{i0}$.}
Let $\Gamma \in G_d^{i0}$. 
Let $v\in \text{Vert}(\Gamma)$ be the vertex at which the
$2^{nd}$ marking  is incident. $D(v)$ is collapsed to $p_i$.
The restriction of $\psi_2$ to
$\overline{M}_\Gamma$ carries the trivial $\mathbf{T}$-action.
Hence, a simple nilpotency result holds: 
\begin{equation}
\label{rrrrrr}
\psi^{\text{dim}(v)+1}_2=0 \in H^*_\mathbf{T}(\overline{M}_\Gamma)
\end{equation}
where $\text{dim}(v)= 
\text{val}(v)-3$ is the dimension of $\overline{M}_{0, \text{val}(v)}$.
The valence bound $\text{val}(v)\leq d+2$ is easily deduced for $\Gamma\in
G_d^{i0}$. It is achieved for graphs with $d$ edges and both
markings incident to $v$. The order of nilpotency of $\psi_2$
is thus bounded
by $d$.  
\vspace{+12pt}

\noindent{\bf{Type} $G_d^{i1}$.}
Let $\Gamma \in G^{i1}_d$. Again
let $v\in \text{Vert}(\Gamma)$ be the vertex at which the
$2^{nd}$ marking is incident. In this case, $D(v)$ is a point.
The vertex 
$v$ is incident to a unique edge $e$ of $\Gamma$.  
Let $e$ connect vertices $v$, $v'$.
Let $j=\mu(v')$.
If $\delta(e) <d$,
let $\Gamma_j$ be the
$2$-pointed graph obtained by contracting $e$:
$\Gamma_j$ is the complete subgraph of $\Gamma$ not containing
$v$ with the $2^{nd}$ marking placed at $v'$.
The graph $\Gamma_j$ is an element of $G^j_{d-\delta(e)}$.
Note also that $|\text{Aut}(\Gamma)|= |\text{Aut}(\Gamma_j)|$.
It is this pruning of graphs which provides recursion
relations for the correlators $z_i$ from equivariant
localization.

The first application of localization is the following
basic  observation:
\begin{lm}
\label{yesy}
The correlators $Z_i(q, \hb, l\leq m+1)$ are
naturally elements of the ring $\Q(\lambda,  \hb)[[q]]$:
$$Z_i(q,\hb,l)= 1 + \sum_d q^d \ \zeta_{id}(\lambda, \hb).$$
Moreover, the rational functions $\zeta_{id}(\lambda, \hb)$
are regular at all the values 
$$\hb= \frac{\lambda_i - \lambda_j}{n}$$ where
 $i\neq j$ and $n\geq 1$.
\end{lm}
\bpf A priori, $\zeta_{id}(\lambda, \hb) \in \R[[\hb^{-1}]]$:
\begin{equation}
\label{zetdef}
\zeta_{id}(\lambda,\hb)= \sum_{k=0}^\infty 
\int_{\overline{M}_d} 
\frac{\psi_2^k E'_{d,2} }{\hb^{k+1}}e_2^*(\phi_i).
\end{equation}
Let $\Gamma \in G_d$. The contribution of $\Gamma$ to
$\zeta_{id}$ is :
$$\text{Cont}_{\Gamma}(\zeta_{id})=\sum_{k=0}^\infty 
\int_{\overline{M}_\Gamma} 
\frac{\psi_2^k E'_{d,2}}{\hb^{k+1} N_\Gamma}e_2^*(\phi_i).$$
By the Type $G_d^{i*}$  vanishing, we obtain:
\begin{equation}
\label{ththth}
\zeta_{id} = \sum_{\Gamma\in G_d^{i0}} \text{Cont}_\Gamma(\zeta_{id})
+ \sum_{\Gamma\in G_d^{i1}} \text{Cont}_\Gamma(\zeta_{id}).
\end{equation}
Let $\Gamma\in G_d^{i0}$. By the Type $G_d^{i0}$ nilpotency
condition, we see:
\begin{equation}
\label{grgrgr}
\text{Cont}_\Gamma(\zeta_{id})=\sum_{k=0}^{d-1} 
\frac{P_{\Gamma,k}(\lambda)}{\hb^{k+1}}
\end{equation}
where $P_{\Gamma,k}(\lambda) \in \Q(\lambda)$.
Let $\Gamma \in G^{i1}_d$. The restriction of $\psi_2$ to
$\overline{M}_\Gamma$ is topologically trivial
with $(\lambda_j - \lambda_i)/\delta(e)$ as equivariant class
(we adhere to the notation of Type $G^{i1}_d$ above). 
Hence, the contributions of $\Gamma$ to
the terms $k\geq 0$ form a geometric series. The series sum is:
\begin{equation}
\label{lucky}
\text{Cont}_\Gamma(\zeta_{id}) =
\frac{P_\Gamma(\lambda)} {(\hb+ \frac{\lambda_i-\lambda_j}{
\delta(e)})},
\end{equation}
where $P_\Gamma(\lambda)\in \Q(\lambda)$.  By 
equations (\ref{ththth}-\ref{lucky}),
$\zeta_{id}\in \Q(\lambda,\hb)$.
The explicit forms of (\ref{grgrgr}) and
(\ref{lucky}) prove the regularity
claim at $\hb=(\lambda_i-\lambda_j)/n$. \epf

The contributions of $G^{i0}$ and 
$G^{i1}$ to the integrals in (\ref{hobo}) yield linear recursion relations
for the correlators $z_i(Q,\hb, l\leq m)$. The contribution of $G^{i0}$
will be the initial part of the relation.
This contribution is analyzed first.

\begin{lm} The contribution $C_i(Q,\hb, l)$ of graph type $G^{i0}$  to
$z_i(Q, \hb,l)$ is determined in cases (i) and (ii) by the following
results:
\begin{eqnarray}
C_i(Q, \hb, l<m) & = & 0 \\
C_i(Q, \hb, m) & = &
-1+ {\text{\em exp}} 
\big( -m!Q + \frac{(m\lambda_i)^m}{\Pi_{\alpha\neq i}(\lambda_i
- \lambda_\alpha)}Q \big), 
\end{eqnarray}
\end{lm}

\bpf
Let $\Gamma \in G_d^{i0}$. 
We follow the notation of Type $G_d^{i0}$ above.
For $d>0$ and $l<m$, we see $(m+1-l)d-1 \geq d$. 
Hence, the restriction of the integrand of (\ref{hobo}) to $\overline
{M}_{\Gamma}$
vanishes by (\ref{rrrrrr}) and the valence bound in these cases.
Thus,  $C_i(Q, \hb, l<m)=0$.

In case $l=m$, the $(1/\hb)$-expansion of the integrand  (\ref{hobo})
contains only one possibly nonvanishing term after restriction to
$\overline{M}_{\Gamma}$:
$\psi^{d-1}_2 E_{d,2}'e_2 ^*(\phi_i)$.
This term also vanishes unless the valence of $v$ is $d+2$.
As previously remarked,
the graphs $\Gamma\in G_{d}$ with valence $d+2$ at a vertex $v$ 
are particularly simple:
they must have $d$ edges (all of degree 1) and both markings incident at
the vertex $v$. It is then straightforward
to explicitly compute the contribution
$C_i(Q, \hb, m)$ from the combinatorics of these simple
graphs via the localization formula. This is the first localization
computation needed by Givental. \epf

The linear recursion relations for $z_i(Q, \hb,l\leq m)$ are
given by:
\begin{equation}
\label{bob}
z_i(Q,\hb,l\leq m)  = 1 + C_i(Q,\hb,l) + \sum_{j\neq i}\sum_{d>0}
Q^{d} \ C_i^j(d, \hb, l) \ z_j(Q, \frac{\lambda_j -\lambda_i}{d}, l),
\end{equation}
where the recursion coefficients are:
\begin{equation}
\label{blob}
C_i^{j}(d, \hb, l) = \frac{1}{\frac{\lambda_i-\lambda_j}{\hb}+d}
\ \frac{\Pi_{r=1}^{ld} \ \frac{ld\lambda_i}{\lambda_j-\lambda_i}+ r}
{\Pi_{\alpha=0}^m \Pi_{r=1,}^{d} {}_{(\alpha,r) \neq (j,d)}
\ \frac{d(\lambda_i- \lambda_\alpha)}{\lambda_j-\lambda_i}+r}.
\end{equation}
The initial term $C_i$ is the contribution of $G^{i0}$ to $z_i$.
The double sum is the contribution of $G^{i1}$. 
A truly remarkable
feature of (\ref{bob}) is the $\hb= (\lambda_j-\lambda_i)/d$
substitution on the right. 
This substitution is well defined by Lemma \ref{yesy}.
Its origin is a normal
bundle factor in the localization formula. 
Let $\Gamma \in G^{i1}_d$. We follow the notation of
Type $G^{i1}_d$ above.  If $\delta(e)=d$, then 
the contribution of $\Gamma$ to $z_i$ equals $Q^d C_i^j(d, \hb,l)$.
Assume $\delta(e)<d$. 
Let $\Gamma_j$ be
the contracted graph obtained from $\Gamma\in G^{i1}_d$ 
as described in Type $G^{i1}_d$ above.
The linear recursion is obtained by the following
equation:
\begin{equation}
\label{contty}
\text{Cont}_{\Gamma}(z_i(Q,\hb,l\leq m)) =Q^{\delta(e)}
C_i^j(\delta(e), \hbar, l) \cdot  \text{Cont}_{\Gamma_j}(
z_j(Q, \frac{\lambda_j-\lambda_i}{\delta(e)}, l)),
\end{equation}
where $\text{Cont}_{\Gamma}$ denote the contribution of
$\Gamma$ to the argument.
Equations (\ref{bob} -\ref{blob}) are deduced directly 
from (\ref{contty}) by summing over all graphs
$\Gamma \in G^{i1}$.
The flag $(v', e)$ in the graph $\Gamma$ corresponds to
a node in the domain curve.  
The normal bundle of $\overline{M}_\Gamma \subset \overline{M}_d$
has a line bundle quotient obtained from the
deformation space of this node. 
This nodal deformation is absent in the normal bundle contributions
for the graph $\Gamma_j$, but appears algebraically in
the evaluation of the
correlator $z_j$ at $\hb= (\lambda_j-\lambda_i)/d$.
Once this graph pruning strategy is noticed and the
explicit recursions given, it is nothing more than
an algebraic check to prove equation (\ref{contty}) from
the localization formulas.  This is the second needed
localization computation.

Define the correlators $Z^*_i\in \R[[\hb^{-1},q]]$ by
\begin{equation}
\label{wqwq}
Z_i^*(q,\hb,l\leq m+1)= \sum_{d=0}^\infty q^d\ 
\frac{\Pi_{r=1}^{ld}  (l
\lambda_i + r\hb)}{\Pi_{\alpha=0}^m 
\Pi_{r=1}^d  (\lambda_i- \lambda_\alpha+r\hb)}.
\end{equation}
For all $l\leq m+1$, 
$Z_i^* \in \Q(\lambda, \hb)[[q]]$, and
the correlators $Z_i^*$
satisfy the regularity property of Lemma \ref{yesy}.
Let $z_i^*(Q, \hb, l\leq m)= Z_i^*(Q\hb^{(m+1-l)}, \hb, l)$. 
A direct algebraic computation shows the
correlators
$z_i^*(Q, \hb, l<m)$ satisfy the recursions (\ref{bob}-\ref{blob}). 
Hence,
$z_i^*(Q,\hb,l<m)= z_i(Q, \hb, l<m)$ since the
recursions clearly have a unique solution.

Define the correlators
$S^*_{\mathbf{T}}\in H_{\mathbf{T}}^*(\proj^m)[[\hb^{-1},t,e^t]]$ by:
\begin{equation}
\label{qertt}
S^*_{\mathbf{T}}(t,\hb,l\leq m+1)= \sum_{d\geq 0} \frac{ e^{(H/\hb+d)t} 
\Pi_{r=0}^{ld}( lH+r\hb)}{
\Pi_{\alpha=0}^{m}
\Pi_{r=1}^d(H-\lambda_\alpha+r\hb)}.
\end{equation} 
The following equality holds:
\begin{equation}
\label{ertt}
\langle \phi_i , S^*_{\mathbf{T}}(t, \hb,l\leq m+1)\rangle 
=e^{\lambda_i t/\hb} l\lambda_i Z_i^*(e^t, \hb, l).
\end{equation}
For $l<m$, we have in addition:
\begin{eqnarray*}
e^{\lambda_i t/\hb} l\lambda_i Z_i^*(e^t, \hb, l<m) 
& = & e^{\lambda_i t/\hb} l\lambda_i Z_i(e^t, \hb,l) \\ 
& = & \langle \phi_i , S_{\mathbf{T}}(t, \hb,l)\rangle.
\end{eqnarray*}
Hence, by property (\ref{vbvb}),
\begin{equation}
\label{anz}
S^*_{\mathbf{T}}(t, \hb, l<m)= S_{\mathbf{T}}(t,\hb,l).
\end{equation}
The non-equivariant correlator $S_X$ is defined by:
$$S_X= \frac{1}{lH} \ S_{\mathbf{T}}(t,\hb,l\leq m+1)
|_{\lambda_i=0, \ \hb=1}$$
(see Section \ref{vbvbvb}). 
Case (i) of the main result (Section 0) is then proven
by (\ref{anz}) and a calculation of the non-equivariant
restriction of $S^*_{\mathbf{T}}$.

For $l=m$, a direct calculation shows the slightly modified
correlator $e^{-m!Q}z_i^*(Q, \hb,m)$ satisfies (\ref{bob}-\ref{blob}). 
Therefore,
$e^{-m!Q}z_i^*(Q, \hb,m)=z_i(Q, \hb,m)$ by uniqueness. 
The equality $e^{-m!e^t/\hb}S^*_{\mathbf{T}}(t,\hb,m)= S_{\mathbf{T}}(t,\hb,m)$
then follows analogously. Givental's relationship (ii) is 
obtained.

By Lemma \ref{wwww}, relations in the quantum cohomology ring
$QH^*_s(X)$ of the hypersurface $X$ are obtained
from differential operators annihilating $S(t,\hb, l)$.
The differential equation:
$$\Big(\hb \frac{d}{dt} \Big)^{m} \ S(t, \hb, l<m) =
 l e^t \Pi_{r=1}^{l-1} \Big( l \hb\frac{d}{dt} +r\hb\Big)\  S(t, \hb,l)$$
yields the relation: $H^m = l^l q H^{l-1}$
in $QH^*_s(X_{l<m})$. Restricted cases of this relation
where proven previously by Beauville [Be].
Similarly, the equation:
$$\Big(\hb \frac{d}{dt}+m! e^t \Big)^{m} S(t, \hb, m) =
m e^t \Pi_{r=1}^{m-1}  
\Big( m \hb \frac{d}{dt}+m m!e^t +r\hb\Big) \ S(t, \hb,m)$$
yields the relation: $(H+m!q)^m = m^m q (H+m!q)^{m-1}$ in
$QH^*_s(X_{m})$.

\section{\bf Case (iii): Calabi-Yau}
\subsection{Linear recursions}
For the Calabi-Yau case $l=m+1$,
no initial terms vanish in (\ref{zdef}) for dimensional reasons. 
It will
be useful to write $Z_i$ in a partially expanded form. 
\begin{eqnarray}
\label{sidd}
Z_i(q, \hb, m+1) &  = &  1 + \sum_{d>0} q^d \Big( \sum_{k=0}^{d-1}
\hb^{-k-1} \int_{\overline{M}_d} \psi_2^k E'_{d,2} e_2^*(\phi_i) \Big) \\
& & \label{fred}
\ \  +\sum_{d>0} 
\big(\frac{q}{\hb}\big)^d 
\int_{\overline{M}_d} 
\frac{\psi_2^{d} }{\hb-\psi_2} E'_{d,2} e_2^*(\phi_i)
\end{eqnarray}
In this case, the definition $z_i(Q, \hb, m+1)= Z_i(Q\hb, \hb, m+1)$ is made.
The treatment of the $\hb$ factor differs in case (iii) 
from the treatment given in (i) and (ii).

The linear recursions for $z_i$ in the Calabi-Yau case take the form:
\begin{eqnarray*}
z_i(Q,\hb,m+1) & = & 1 + \sum_{d>0} \frac{Q^d}{d!} R_{id} \\ 
& & \ \ +
\sum_{d>0}\sum_{j\neq i} 
Q^d \ C_i ^j(d,\hb,m+1) \ z_j(Q,\frac{\lambda_j-\lambda_i}
{d}, m+1),
\end{eqnarray*}
where $R_{id}= \sum_{j=0}^{d} R_{id}^j \hb^{d-j}$ is a
polynomial in $\R[\hb]$ of $\hb$-degree at most $d$, and
the recursion coefficient $C_i^j$ is determined by:
\begin{equation}
\label{cycof}
C^j_i(d,\hb,m+1)= \frac{1}{\lambda_i-\lambda_j+d\hb} \
\frac{\Pi_{r=1}^{(m+1)d} (m+1)\lambda_i + r\frac{\lambda_j-\lambda_i}{d}}
{d! \Pi_{\alpha\neq i} \Pi_{r=1, \ (\alpha,r)\neq (j,d)}^d \ 
\lambda_i-\lambda_{\alpha}+ r\frac{\lambda_j-\lambda_i}{d}}.
\end{equation} 
The proof this recursion relation is similar (but not
identical) to the proofs in
cases (i) and (ii) of (\ref{bob}).
Equation (\ref{bob}) was derived by separating the
contribution of graph types $G^{i0}$ and $G^{i1}$.
Here, we instead separate the contributions of the 
terms (\ref{sidd}) and (\ref{fred}) in the expansion
of $Z_i(q, \hb,m+1)$. 
The contribution of (\ref{sidd}) is easily seen to be of the
form $1 + \sum_{d>0} Q^d R_{id}/ d!$
given by the initial term in the recursion -- the polynomials
$R_{id}$ are not specified. 
Next, the contribution of (\ref{fred}) is analyzed.
The 
graph type $G^{i0}$
contributions to (\ref{fred}) vanish by the argument in Section 3 since
$\psi_2$-degree
is too high. Hence, only graphs of type $G^{i1}$  contribute to (\ref{fred}).
The graph pruning strategy is now applied as before. 

In this case,
we include the details of the required localization
calculation.  
Let $\Gamma \in G_d^{i1}$.
Two equations are needed.
Let $\Gamma$ be the unique graph of type $G_d^{i1}$
with a single edge $e$ connecting fixed points $i$ and $j\neq i$
and satisfying $\delta(e)=d$. The first equation is:
\begin{equation}
\label{ffrr}
\text{Cont}_{\Gamma} \Big( Q^d \int _{\overline{M}_d}
\frac{\psi_2^{d}}{\hb-\psi_2} E'_{d,2} e_2^*(\phi_i)\Big)
= Q^d C_i^j(d,\hb, m+1).
\end{equation}
The proof is by a computation of the  left  contribution.
The stack $\overline{M}_\Gamma$ is 0 dimensional with
$|\mathbf{G}|=d$; the space $\widetilde{M}_\Gamma$ is a
regular point. 
Let $$\mu: (C, x_1, x_2) \rarr \proj^m$$ be the fixed map 
corresponding to $\Gamma$.
The restriction of the 
equivariant top Chern class (or Euler class)
 of $E'_{d,2}$  to $\overline{M}_\Gamma$ is:
$$c_{\text{top}}(E'_{d,2})|_{\overline{M}_\Gamma}
= \Pi_{r=1}^{(m+1)d} (m+1) 
\lambda_i + r\frac{\lambda_j- \lambda_i}{d}.$$
The Euler class of $N_\Gamma$ is obtained from 
$H^0(C, \mu^*(T_{\proj^m}))$
(after subtracting the trivial weight obtained from the
unique infinitesimal automorphism  fixing $x_1$ and $x_2$).
The weights of $H^0(C, \mu^*(T_{\proj^m}))$ are determined by the
$\mu$ pull-back of the Euler sequence:
\begin{eqnarray*}
c_{\text{top}}(N_\Gamma) & = & \Pi_{\alpha=0}^m \Pi_{r=0, \ 
(\alpha,r)\neq (i,0), (j,d)}^d
\ \lambda_i- \lambda_\alpha + r \frac{\lambda_j-\lambda_i}{d} \\
& = & d! (\frac{\lambda_j-\lambda_i}{d})^d \cdot 
\Pi_{\alpha\neq i} (\lambda_i-\lambda_\alpha) \cdot \\& & 
\Pi_{\alpha\neq i} \Pi_{r=1, \ 
(\alpha,r)\neq (j,d)}^d
\ \lambda_i- \lambda_\alpha + r \frac{\lambda_j-\lambda_i}{d}.
\end{eqnarray*}
The  classes $\psi_2$ and $e_2^*(\phi_i)$ restrict to
$(\lambda_j-\lambda_i)/d$ and $\Pi_{\alpha\neq i}(\lambda_i 
-\lambda_\alpha)$ respectively. Since,
$$
\text{Cont}_{\Gamma} \Big( Q^d \int _{\overline{M}_d}
\frac{\psi_2^{d}}{\hb-\psi_2} E'_{d,2} e_2^*(\phi_i)\Big)
= \frac{Q^d}{|\mathbf{G}|} \int _{\widetilde{M}_\Gamma}
\frac{\psi_2^{d}}{(\hb-\psi_2) \cdot  N_\Gamma} E'_{d,2} e_2^*(\phi_i),
$$
equation (\ref{ffrr}) is
an algebraic  consequence of these weight calculations
(pulled back to $H^*_{\mathbf{T}} ({\widetilde{M}_\Gamma})$).
 
Next, let $\Gamma \in G^{i1}_d$ satisfy $\delta(e)<d$
(following the notation of Type $G^{i1}$ in Section 3).
The second equation is:
\begin{equation}
\label{sscc}
\text{Cont}_{\Gamma} \Big( Q^d \int _{\overline{M}_d}
\frac{\psi_2^{d}}{\hb-\psi_2} E'_{d,2} e_2^*(\phi_i)\Big) =
\end{equation}
$$ Q^{\delta(e)} C_i^j(\delta(e), \hb,m+1) \cdot
\text{Cont}_{\Gamma_j} \Big( 
(\frac{\lambda_j-\lambda_i}{\delta(e)} Q)^{d-\delta(e)}
\int _{\overline{M}_{d-\delta(e)}}
\frac{  E'_{d-\delta(e),2}  }{\frac{\lambda_j-\lambda_i}{\delta(e)}
-\psi_2}  
e_2^*(\phi_j)\Big).
$$
The proof is by expanding the left contribution.

Standard weight calculations (via a natural restriction
sequence of sections of $\mu^*(\oh_{\proj^m}(m+1))$ 
to  the  component $D(e)$) yield the following product
formula of weights:
\begin{equation}
\label{qkqk}
c_{\text{top}} (E'_{d,2})|_{\overline{M}_\Gamma}= 
c_{\text{top}} (E'_{d-\delta(e),2})|_{\overline{M}_{\Gamma_j}}
\cdot  \Pi_{r=1}^{(m+1)\delta(e)} (m+1) 
\lambda_i + r\frac{\lambda_j- \lambda_i}{\delta(e)}.
\end{equation}

The normal bundle $N_\Gamma$  is determined in K-theory
as a sum of two pieces. 
Let $[\mu:(C, x_1, x_2) \rarr \proj^m] \in \overline{M}_\Gamma$.
The first piece is topologically trivial
with weight obtained from the representation
$H^0(C, \mu^*(T_{\proj^m})) $ after removing the
infinitesimal automorphisms. The second piece 
is a direct
sum of line bundles 
obtained from the deformation spaces of  the nodes
of $C$ forced by $\Gamma$. Recall, the
second marking of $\Gamma_j$ corresponds to
the forced node of $C$ lying on $D(e)$.

The normal bundle piece obtained from $H^0(C, \mu^*(T_{\proj^m}))$
may be decomposed via restriction to $D(e)$ as in (\ref{qkqk}). 
Note that $\widetilde{M}_\Gamma$ and $\widetilde{M}_{\Gamma_j}$
are canonically isomorphic. Via this isomorphism, we find 
\begin{equation}
\frac{e_2^*(\phi_i)}{N_\Gamma} =\frac{1}{\frac{\lambda_j-\lambda_i}
{\delta(e)}
-\psi_2} \frac{e_2^*(\phi_j)}
{N_{\Gamma_j}}  \cdot
\end{equation}
$$
\frac{1}{\delta(e)! (\frac{\lambda_j-\lambda_i}{\delta(e)})^{\delta(e)}
\cdot \Pi_{\alpha\neq i} \Pi_{r=1, \ 
(\alpha,r)\neq (j,\delta(e))}^{\delta(e)}
\lambda_i- \lambda_\alpha + r \frac{\lambda_j-\lambda_i}{\delta(e)}},
$$ 
where the left and right sides 
are a naturally classes on $\widetilde{M}_\Gamma$
and $\widetilde{M}_{\Gamma_j}$ respectively.
The first term on the right is the nodal deformation
corresponding to the pruned node.

It is important to realize the treatment of the second marking
differs for $\widetilde{M}_\Gamma$ and $\widetilde{M}_{\Gamma_j}$.
The natural pull-back of $\psi_2$ to $\widetilde{M}_\Gamma$
is of pure weight $(\lambda_j-\lambda_i)/\delta(e)$. 

Finally, we have $|\mathbf{G}_\Gamma|= \delta(e) |\mathbf{G}_{\Gamma_j}|$.
As the $\Gamma$ and $\Gamma_j$ contributions in (\ref{sscc})
may be integrated on the tilde space  (with automorphism
corrections),  equation (\ref{sscc}) now follows algebraically.
The linear recursions are obtained from (\ref{ffrr}) and
(\ref{sscc}) by summing over graphs of type $G^{i1}$.
This is the third use of the full localization formulas for
the moduli space of maps.

\subsection{Polynomiality}
\label{shc}
The Calabi-Yau case is difficult for several reasons.
The recursion relations for $z_i$ are not yet 
determined as the functions $R_{id}$ are unknown. 
It is necessary to find additional conditions satisfied 
by the correlators $z_i$. Givental's idea here is to
prove a 
polynomiality constraint satisfied by a related 
double correlator $\Phi$. Define $\Phi(z,q)\in 
\Q(\lambda,\hb)[[z,q]]$ 
by:
\begin{equation}
\label{zdoub}
\Phi(z,q)= \sum_{i=0}^m \frac{(m+1)\lambda_i}{\Pi_{j\neq i}
(\lambda_i- \lambda_j)}
\ e^{\lambda_iz} Z_i(qe^{z\hb}, \hb,m+1) Z_i(q, -\hb,m+1).
\end{equation}
A constraint on $\Phi(z,q)$ may be interpreted
as a further condition on the correlators $z_i$.

A geometric construction is needed for the polynomiality
constraint. Consider a new $1$-dimensional torus $\com^*$.
Let $\mathbb{Q}[\hb]$ be the standard presentation
of the equivariant cohomology ring of $\com^*$
(again, $\hb$ is the first Chern class of the dual of the
standard representation of $\com^*$). Let $\com^*$
act on the vector space $V=\com^2$ via the exponential
weights $(0,-1)$. Let $y_1, y_2$ be the respective
fixed points for the induced action on $\proj^1=\proj(V)$. 
The equivariant Chern classes of the tangent representations
at the fixed points are
$\hb, -\hb$ respectively. 
Recall from Section \ref{aaction} the $\mathbf{T}$-action
on $W$.
There are naturally
induced $(\com^* \times \mathbf{T})$-actions on
$\proj(V) \times \proj(W)$ and $\overline{M}_{0,2}(\proj(V)
\times \proj(W), (1,d))$.
The space of interest to us will be: 
$$L_d = e_1^{-1}\big(\{y_1\} \times \proj(W) \big) 
\cap e_2^{-1}\big(\{y_2\} \times \proj(W)\big) \subset
\overline{M}_{0,2}(\proj(V)
\times \proj(W), (1,d)).$$
$L_d$ is easily seen to be a nonsingular, $(\com^* \times 
\mathbf{T})$-equivariant substack.

Let $L'_d$ denote the
polynomial space $\proj(W \otimes Sym^d(V^*))$ with the
canonical $(\com^* \times \mathbf{T})$-representation. 
A degree $d$ algebraic map $\proj(V) \rarr \proj(W)$
canonically yields a point in $L_d'$. 
There
is a natural $(\com^*\times \mathbf{T})$-equivariant
morphism
$$\mu:M_{0,2}(\proj(V)\times \proj(W), (1,d))
\rarr L'_d$$
obtained by identifying an element of the left
moduli space with the graph of a uniquely
determined map $\proj(V) \rarr \proj(W)$.
It may be shown that $\mu$ extends to
a $\com^*\times \mathbf{T}$-equivariant morphism from the stack
$\overline{M}_{0,2}(\proj(V)\times \proj(W), (1,d))$ [G1], [LLY].
Let $\mu: L_d \rarr L'_d$ be the induced map.
Let $P\in H^*_{\com^*\times \mathbf{T}}(L_d')$ be the
first Chern class of $\oh_{L_d'} (1)$.
Let $E_d$
be the equivariant bundle on $L_d$ with fiber over a
stable map $[(f_V\times f_W) :C \rarr P(V) \times P(W)]$ equal
to $H^0(C,f_W^*(\oh_{\proj(W)}(m+1)))$.

\begin{lm}
\label{zaza}
There is an equality:
\begin{equation}
\label{rdrd}
\Phi(z,q)= \sum_{d\geq 0} q^{d} \int_{L_d}
e^{\mu^*(P)\cdot z} E_d,
\end{equation}
where the integral on the right is the $(\com^*\times
\mathbf{T})$-equivariant push forward to a point. 
\end{lm}

\bpf
This is the fourth (and last) localization calculation on the
moduli space of maps
needed by Givental. 
The remarkable feature of this equality is the following.
On the left side of (\ref{rdrd}),
$\hb$ is a formal parameter. On the right side, it
an element of equivariant cohomology.
As $L_d$ is a nonsingular
stack, the $\com^*\times \mathbf{T}$-localization formula
yield an explicit graph summation answer for the
integral on the right which
is directly matched with (\ref{zdoub}).

The first step is identify the graph types of the
fixed loci of $L_d$. 
Recall the definitions of $G^{i0}$ and $G^{i1}$ from Section 3.
Let $G^i= G^{i0} \scup G^{i1} \scup \{\text{Triv}(i)\}$
where $\text{Triv}(i)$ is the edgeless two pointed
graph with  a single vertex $v$ satisfying $\mu(v)=p_i$.
Let $\text{deg}(\text{Triv(i)})=0$.
The components of $L_d^{\com^*\times \mathbf{T}}$ 
are
in bijective correspondence to triples $(i,\Gamma_1, \Gamma_2)$
where $0\leq i \leq m$  and $\Gamma_1 , \Gamma_2\in G^{i}$  
satisfy $\text{deg}(\Gamma_1)+\text{deg}(\Gamma_2)=d$.
The graphs $\Gamma_1, \Gamma_2$ describe the configurations
lying over the points $y_1, y_2 \in \proj(V)$ respectively.
A fixed map $$\mu:(C, x_1,x_2) \rarr \proj(V)\times \proj(W)$$
in the corresponding component satisfies the following
properties. The domain is a union of three subcurves
$C=C_1\scup C_m \scup C_2$.  The curve $C_m$ is mapped
isomorphically by $\mu$ to $\proj(V) \times \{p_i\}$.
$C_1$ and $C_2$ contain $x_1$ and $x_2$ 
and lie over $y_1$ and $y_2$ respectively.
The Lemma will follow from the calculation of
the contribution of $(i, \Gamma_1, \Gamma_2)$ to the 
integral in (\ref{rdrd}).

Let $\mathbf{\Gamma}=(i, \Gamma_1, \Gamma_2)$.
Let $d_1, d_2$ equal $\text{deg}(\Gamma_1),
\text{deg}(\Gamma_2)$ respectively.
We treat the generic case:  $d_1, d_2 >0$. 
The degenerate cases in which either $\Gamma_1$ or $\Gamma_2$ equals
$\text{Triv}(i)$ are computed analogously.
The contribution equation is:
\begin{eqnarray*}\label{helel}
\text{Cont}_{\mathbf{\Gamma}}( q^d \int_{L_d} e^{\mu^*(P) \cdot z}
E_d) &= & \frac{(m+1) \lambda_i}{\Pi_{\alpha\neq i} 
\lambda_i-\lambda_\alpha} e^{\lambda_iz}\cdot \\
& &  (qe^{z\hb})^{d_1} \cdot
\text{Cont}_{\Gamma_1}( \int_{\overline{M}_{d_1}}
\frac{E'_{{d_1},2}}{\hb-\psi_2} e_2^*(\phi_i) ) \cdot \\
& & q^{d_2} \cdot  \text{Cont}_{\Gamma_2}( \int_{\overline{M}_{d_2}}
\frac{E'_{{d_2},2}}{-\hb-\psi_2} e_2^*(\phi_i) ).
\end{eqnarray*} 
The contribution equation in the degenerate cases is
identical (with the convention $\text{Cont}_{\text{Triv(i)}}=1$).

The equation is proven by expanding the localization formula
for the left contribution.  Note first that the fixed stack
$\overline{M}_{\mathbf{\Gamma}} \subset L_d$ is
naturally isomorphic to $ \overline{M}_{\Gamma_1}
\times \overline{M}_{\Gamma_2}$. As 
$$\mu(\overline{M}_{\mathbf{\Gamma}})= [\com_i \otimes 
[(y_1^*)^{d_2} (y_2^*)^{d_1}]],$$ the class
$\mu^*(P)$ is pure weight equal to $\lambda_i + d_1 \hb$.
The class $c_{\text{top}}(E_d)|_{\overline{M}_{\mathbf{\Gamma}}}$
is pure weight and factors as:
$$(m+1)\lambda_i \cdot c_{\text{top}} (E'_{d_1,2})
|_{\overline{M}_{\Gamma_1}} \cdot
c_{\text{top}}(E'_{d_2,2})|_{\overline{M}_{\Gamma_2}}$$
by the restriction sequence to $C_m$. Similarly,
$$\frac{\Pi_{\alpha\neq i} \lambda_i-\lambda_\alpha}
{N_{\mathbf{\Gamma}}}$$ is computed to equal
the product of $e_2^*(\phi_i)/((\hb-\psi_2) N_{\Gamma_1})$
from $\overline{M}_{\Gamma_1}$ with
  $e_2^*(\phi_i)/((-\hb-\psi_2) N_{\Gamma_2})$ from
$\overline{M}_{\Gamma_2}$. This normal bundle expression
is obtained by the restriction sequence of tangent
sections to $C_m$ and an accounting of nodal deformations.
As the $N_{\mathbf{\Gamma}}$ is the normal bundle in $L_d$,
only tangent sections of $H^0(C,  \mu^*(\proj(V))$ vanishing
at the markings $x_1$ and $x_2$ appear in the normal bundle
expression.  
The contribution equation now follows directly.

Equation (\ref{rdrd}) is obtained from the contribution
equation, the definition of $Z_i(q,\hb,m+1)$, and
a sum over graphs.  
\epf

By Lemma \ref{zaza}, $\Phi(z, q)$ may be rewritten as:
\begin{equation}
\label{ppolyy}
\Phi(z, q)= 
\sum_{d\geq0} q^{d}
\int_{L'_d} e^{P z} \mu_*(c_{\text{top}}(E_d)).
\end{equation}
The group
$\com^*\times \mathbf{T}$ acts with $(m+1)(d+1)$
isolated fixed points on $L_d'$. 
A weight calculation of the representation
$W\otimes Sym^d(V^*)$ yields the standard 
presentation:
$$H^{*}_{\com^*\times\mathbf{T}}(L_d')=
\mathbb{Q}[P,\lambda, \hb]/ (\Pi_{\alpha=0}^m \Pi_{r=0}^d
(P - \lambda_\alpha -r \hb)).$$
As
 $\mu_*(c_{\text{top}}(E_d)) \in H^{(m+1)d+1}_{\com^*\times
\mathbf{T}}(L_d')$,
there is a unique polynomial $$E^{Z}_d(P,
\hb,\lambda)\in \Q[P,\lambda, \hb]$$ 
of homogeneous degree $(m+1)d+1$ satisfying
$\mu_*(c_{\text{top}}(E_d))=E^{Z}_d(P,\lambda,\hb)$ in
$H^*_{\com^*\times
\mathbf{T}}(L_d')$.
The Bott residue formula for the  integral in (\ref{ppolyy})
then yields:
\begin{equation}
\label{polycon1}
\Phi(z, q)= 
\frac{1}{2\pi i} \oint  e^{P z} \sum_{d\geq 0}
\frac{ q^{d} E_d^{Z}(P, \lambda,\hb)}
{ \Pi_{\alpha=0}^m \Pi_{r=0}^d
(P - \lambda_\alpha -r \hb)} dP.
\end{equation}
Givental's polynomiality constraint is the following:
$\Phi(z, q)$ is expressible as a
residue integral of the form (\ref{polycon1})
where $E_d^{Z}(P, \lambda,\hb)\in  
\Q[P,\lambda, \hb]$ is of $P$-degree at most 
$(m+1)d+m$.

\subsection{Correlators of class $\mathcal{P}$}
\label{dede}
Let $\{Y_i(q,\hb)\}_{i=0}^m\subset \R[[\hb^{-1},q]]$ be a set of 
functions (called correlators).
Assume the correlators 
$Y_i$ satisfy the rationality and regularity
conditions of Lemma \ref{yesy}: $Y_i \in \Q(\lambda, \hb)[[q]]$
with no poles at $\hb= (\lambda_i- \lambda_j) /n$ (for
all $j\neq i$ and $n\geq 1$).
Let $y_i(Q,\hb)= Y_i(Q\hb, \hb)$.
Let $y_i$ satisfy the following recursion relation:
\begin{equation}
\label{llll}
y_i(Q,\hb) =  1 + \sum_{d>0} \frac{Q^d}{d!} I_{id}  
+
\sum_{d>0}\sum_{j\neq i} 
Q^d \ C_i ^j(d,\hb,m+1) \ y_j(Q,\frac{\lambda_j-\lambda_i}
{d}),
\end{equation}
where $I_{id}= \sum_{j=0}^{d} I_{id}^j \hb^{d-j}\in \Q(\lambda)[\hb]$ 
is an element of $\hb$-degree at most $d$.
The recursions (\ref{llll}) clearly
determines $y_i$ uniquely from the initial data $I_{id}$.
A direct algebraic consequence of (\ref{llll}) is
the existence of a unique expression:
\begin{equation}
\label{wawa}
y_i(Q,\hb)= \sum_{d\geq 0} Q^d \frac{N_{id}}
{d! \Pi_{j\neq i} \Pi_{r=1}^d (\lambda_i- \lambda_j
+r \hb)},
\end{equation}
where $N_{id}\in \Q(\lambda)[\hb]$ is a polynomial of $\hb$-degree
at most $(m+1)d$, and $N_{i0}=1$.
We may also consider the double correlator $\Phi^Y\in \Q(\lambda,\hb)
[[z,q]]$:
\begin{equation}
\label{mama}
\Phi^Y(z,q)=
\sum_{i=0}^m \frac{(m+1)\lambda_i}{\Pi_{j\neq i}
(\lambda_i- \lambda_j)}
\ e^{\lambda_iz} Y_i(qe^{z\hb}, \hb) Y_i(q, -\hb).
\end{equation}
After the substitution of (\ref{wawa}) in (\ref{mama}),
a straightforward algebraic computation shows:
\begin{equation}
\label{polycon}
\Phi^Y(z, q)= 
\frac{1}{2\pi i} \oint  e^{P z} \sum_{d\geq 0}
\frac{ q^{d} E_d^{Y}(P, \lambda,\hb)}
{ \Pi_{\alpha=0}^m \Pi_{r=0}^d
(P - \lambda_\alpha -r \hb)} dP,
\end{equation}
where 
$E^Y_d= \sum_{k=0}^{(m+1)d+m} f_k(\lambda,\hb) P^k$ is the 
unique function of $P$-degree at most 
$(m+1)d+m$ determined by the values
at the $(m+1)(d+1)$ evaluations $P=\lambda_i+r \hb$ ($0 \leq i
\leq m$, $0 \leq r \leq d$): 
\begin{equation}
\label{hottt}
E^Y_d(\lambda_i+ r \hb)= (m+1) \lambda_i N_{ir}(\hb) N_{i(d-r)}(-\hb).
\end{equation}
In general, the coefficients $f_k(\lambda, \hb)\in \Q(\lambda,\hb)$
will be rational functions.
The correlators $Y_i$ satisfy Givental's polynomiality condition
if $E^Y_d \in \Q[P,\lambda, \hb]$.
\begin{lm}
\label{cker}
The correlators $Y_i$ satisfy Givental's polynomiality condition
if and only if $\Phi^Y(z, q) \in \Q[\lambda,\hb][[z, q]]$.
\end{lm}
\bpf 
By the Bott residue formula, the integral
\begin{equation}
\label{fintegr}
\frac{1}{2\pi i} \oint   \sum_{d\geq 0} 
\frac{ P^k }
{ \Pi_{\alpha=0}^m \Pi_{r=0}^d
(P - \lambda_\alpha -r \hb)} dP
\end{equation}
simply computes the $\com^* \times \mathbf{T}$
equivariant push-forward to a point of the
class $P^k \in H^*_{\com^*\times \mathbf{T}}(L_d')$. 
We therefore see:
\begin{enumerate}
\item[(a)]
for $k< (m+1)d+m$, (\ref{fintegr}) vanishes,
\item[(b)] for 
$k=(m+1)d+m$, (\ref{fintegr}) equals 1,
\item[(c)] 
for $k>(m+1)d+m$, (\ref{fintegr})   is 
an element of $\Q[\lambda,\hb]$. 
\end{enumerate}
Expand the integrand of (\ref{polycon}) 
in power series by $e^{Pz}= \sum_{k=0}^{\infty} (Pz)^k/k!$.
Properties (a)-(c) then prove that the 
polynomiality of the coefficients of $E^Y_d= 
\sum_{k=0}^{(m+1)d+m} f_k(\lambda,\hb) P^k$ is
equivalent to the polynomiality of all coefficients
of the terms $\{z^k q^{d}\}_{k=0}^{\infty}$ 
in $\Phi^Y(z, q)$.\epf

A set of correlators $Y_i\in \R[[\hb^{-1},q]]$ is defined to be of
class $\mathcal{P}$ if the following three conditions are
satisfied.
\begin{enumerate}
\item[I.] The rationality and regularity conditions hold. 
\item[II.] The correlators
$y_i$ satisfy relations of the
form (\ref{llll}).
\item[III.] Givental's polynomiality condition is met. 
\end{enumerate}
A suitable interpretation of II actually implies I, but we
separate these conditions for clarity.

The most important property of class $\mathcal{P}$ is 
Givental's uniqueness result.
\begin{lm}
\label{givu}
Let $Y_i, \overline{Y}_i 
\in \mathbf{R}[[\hb^{-1},q]]$ be two sets of correlators of
class $\mathcal{P}$. 
If
\begin{equation}
\label{assump}
\forall i, \ \ Y_i= \overline{Y}_i \ \ \text{modulo} \ \  \hb^{-2},
\end{equation}
then the sets of correlators
agree identically: $Y_i=\overline{Y}_i$.
\end{lm}

\bpf
Let $I_{id}$ and $\overline{I}_{id}$ be the
respective initial data in the associated recursions (\ref{llll}).
By the recursion formula (\ref{llll}) and the coefficient
formula (\ref{cycof}), we obtain the equality
\begin{equation}
\label{twocof}
Y_i= \sum_{d\geq0} q^d (I^0_{id}+\frac{I^1_{id}}{\hb} )
\ \ \text{modulo} \ \  \hb^{-2}
\end{equation}
(and analogously for $\overline{Y}_i$).
Assumption (\ref{assump}) therefore implies
$I^0_{id}= \overline{I}^0_{id}$ and $I^1_{id}=\overline{I}^1_{id}$
for all $i$ and $d$. In particular, $I_{i1}=\overline{I}_{i1}$.

To establish the Lemma,
it is sufficient to prove $I_{id}=\overline{I}_{id}$ by induction. Assume
$I_{ik}=\overline{I}_{ik}$ for
all $0\leq i\leq m$ and $k<d$. The equality $N_{ik}=\overline{N}_{ik}$
for $k<d$ then follows from the recursions. By (\ref{hottt}), 
$\delta E_d=E_d^{Y}- E_d^{\overline{Y}}$ vanishes at $P= \lambda_i +r \hb$
for all $i$ and $1\leq r \leq d-1$. Hence, the
polynomial $\delta E_d$ is divisible by
$\Pi_{j=0}^m \Pi_{r=1}^{d-1}(P- \lambda_j-r\hb).$
By (\ref{hottt}) and the recursion (\ref{llll}), a computation shows:
$$\delta E_d(P=\lambda_i+d\hb) =(m+1)\lambda_i 
\Pi_{j\neq i} \Pi_{r=1}^d (\lambda_i-\lambda_j+r\hb)
\ (I_{id}- \overline{I}_{id}).$$
By the polynomiality condition $\delta E_d \in \Q[P, \lambda, \hb]$
and the above divisibility, we find $\hb^{d-1}$ divides
$I_{id}- \overline{I}_{id}$. Therefore,
the initial data is allowed to differ only in the
$\hb^d$ and $\hb^{d-1}$ coefficients. However, these
coefficient are precisely the two appearing in (\ref{twocof})
which agree by 
assumption (\ref{assump}). We have proven the equality
$I_{id}=\overline{I}_{id}$. The inductive
step is complete. \epf

By the results of Section \ref{shc},  the correlators
$Z_i(q, \hb, m+1)$ are of class $\mathcal{P}$.
Recall the hypergeometric correlators $Z_i^*(q,\hb, m+1)$ defined
by (\ref{wqwq}).
A straightforward exercise in algebra shows the 
correlators $Z_i^*$ 
also to be of class $\mathcal{P}$. 
The polynomials $E_d^{Z^*}(P, \lambda, \hb)$
associated to the correlators $Z_i^*$ are:
$$E_d^{Z^*}= \Pi_{r=0}^{(m+1)d} ((m+1) P-r \hb).$$
 
The
two sets of correlators $Z_i, Z_i^*$ do not agree modulo $\hb^{-2}$.
The expansions modulo $\hb^{-2}$ may be explicitly evaluated.
From expression (\ref{sidd}), the $\hb^0$ term in
$Z_i$ is 1. The $\hb^{-1}$ term in (\ref{sidd}) vanishes
since the classes in the relevant integrals over 
$\overline{M}_{d}$
are pull-backed via the map forgetting the first marking. Hence,
$Z_i= 1$
modulo $\hb^{-2}$.
A direct computation yields:
$$Z_i^*= F(q)+\frac{ \lambda_i (m+1)(G_{m+1}(q)-G_1(q))+
G_1(q)\sum_{\alpha=0}^m \lambda_\alpha }{\hb}
\ \ \text{modulo} \ \  \hb^{-2},$$
where the functions $F(q)$ and $G_l(q)$ 
are defined by: 
$$F(q)= \sum_{d=0}^\infty q^d \frac{((m+1)d)!}{(d!)^{m+1}},
\ \ 
G_l(q)= \sum_{d=1}^{\infty} q^d \frac{((m+1)d)!}{(d!)^{m+1}}
\Big( \sum_{r=1}^{ld} \frac{1}{r}\Big).$$
The last step in the proof of the Calabi-Yau case (iii)
is the following. An explicit transformation $\overline{Z}_i$
of the  correlator $Z_i$ is found which satisfies:
\begin{enumerate}
\item[(1)] $\overline{Z}_i$ is of class $\mathcal{P}$,
\item[(2)] $\overline{Z}_i = Z^*_i$ modulo $\hb^{-2}$.
\end{enumerate}
Then, by Lemma \ref{givu}, $\overline{Z}_i= Z^*_i$.
This transformation  will yield the
Mirror prediction in the quintic $3$-fold case.

\subsection{Transformations}
Let $Y_i$ be a set of correlators of class $\mathcal{P}$.
We define three transformations:
\begin{enumerate}
\item[(a)] $\overline{Y}_i(q, \hb) = f(q)\ Y_i(q,\hb)$,
\item[(b)] $\overline{Y}_i(q, \hb) = 
\text{exp}(\lambda_i g(q)/\hb)\ Y_i(
q  e^{g(q)}, \hb)$,
\item[(c)] $\overline{Y}_i(q,\hb) = \text{exp}(Cg(q)/\hb)\ Y_i(q,\hb)$,
\end{enumerate}
where $f(q), g(q) \in \Q[[q]]$ satisfy 
$f(0)=1$ and $g(0)=0$, and $C\in \R$ is
a homogeneous linear function of the $\lambda$'s.
\begin{lm}
\label{trann}
In each case (a)-(c), $\overline{Y}_i$ is a 
set of correlators of class $\mathcal{P}$.
\end{lm}

\bpf 
Since rational functions in $\lambda,\hb$ satisfying the
regularity condition of Lemma \ref{yesy} form a subring,
the correlators $\overline{Y}_i$
clearly satisfy condition I of class $\mathcal{P}$.
A direct algebraic check shows the correlators
$\overline{y}_i$ satisfy recursion relations 
of the form (\ref{llll}). The initial terms $\overline{I}_{id}$
change, but remain in $\Q(\lambda)[\hb]$ of $\hb$-degree
at most $d$. The values $f(0)=1$ and $g(0)=0$ are needed
for this verification.
Condition II therefore holds for $\overline{Y}_i$.

Condition III of class $\mathcal{P}$ is checked via
Lemma \ref{cker}. The transformations (a)-(c) have
the following effect on the double correlator:
\begin{enumerate}
\item[(a)] $\Phi^{\overline{Y}}(z,q) = f(qe^{z\hb})
f(q) \cdot \Phi^Y(z ,q)$,
\item[(b)] $\Phi^{\overline{Y}}(z ,q) =
\Phi^Y( z+ (g (qe^{z\hb})-g(q))/\hb , qe^{g(q)})$,
\item[(c)] $\Phi^{\overline{Y}}(z ,q) = \text{exp}(
C\cdot (g(qe^{z\hb})-g(q))/\hb)
 \cdot \Phi^Y(z,q)$.
\end{enumerate}
In each case, $\Phi^{\overline{Y}}$ is easily seen
to remain in $\Q[\lambda,\hb][[z,q]]$.
Case (a) is clear. Since
$$\frac{g(qe^{z\hb})-g(q)}{\hb} 
\in \Q[\lambda,\hb][[z,q]],$$
the change of variables in
case (b)  and multiplication in case (c)
preserve membership in $\Q[\lambda,\hb][[z,q]]$. \epf

The transformation from $Z_i(q, \hb, m+1)$ to
$Z_i^*(q,\hb,m+1)$ can now be established.
Define the correlators $\overline{Z}_i$ by
\begin{equation*}
\overline{Z}_i(q,\hb)=F(q) \cdot \text{exp}(\frac{
(m+1)\lambda_i (G_{m+1}(q)-G_1(q))+ G_1(q)\sum_{\alpha=0}^m 
\lambda_\alpha}{\hb F(q)}) \cdot 
\end{equation*}
$$Z_i(q\cdot \text{exp}(\frac{(m+1)(G_{m+1}(q)-G_1(q))}{F(q)}), \hb, m+1).$$
By a composition of
transformations established in Lemma \ref{trann}, the correlators 
$\overline{Z}_i$  
are of class $\mathcal{P}$. An explicit calculation using
the results of Section \ref{dede}
shows $\overline{Z}_i(q,h) = Z_i^*(q,\hb,m+1)$ modulo $\hb^{-2}$.
By Lemma \ref{givu}, $\overline{Z}_i(q,\hb)=Z_i^*(q, \hb,m+1)$.

Consider the change of variables defined by:
\begin{equation}
\label{nnnnn}
T= t+ \frac{(m+1)(G_{m+1}(e^t)-G_1(e^t))}{F(e^t)}.
\end{equation}
Exponentiating (\ref{nnnnn}) yields
\begin{equation}
\label{nnnnnn}
e^T= e^t\cdot\text{exp}(\frac{(m+1)(G_{m+1}(e^t)-G_1(e^t))}{F(e^t)}).
\end{equation}
Together (\ref{nnnnn}) and (\ref{nnnnnn})
define a change of variables from formal series in $T, e^T$
to formal series in $t, e^t$. This transformation
is easily seen to be invertible.

Let $S_\mathbf{T}(T,\hb,m+1) \in H_\mathbf{T}^*(\proj^m)
[[\hb^{-1},T,e^T]]$ be the equivariant correlator (\ref{lindad})
in the variable $T$. Let the correlator
$\overline{S}_{\mathbf{T}}(t,\hb)\in H_\mathbf{T}^*(\proj^m)
[[\hb^{-1},t,e^t]]$ be obtained from $S_\mathbf{T}(T,\hb,m+1)$
by the change of variables (\ref{nnnnn}) followed by
multiplication by the function
$$F(e^t)\cdot \text{exp}(\frac{G_1(e^t)\sum_{\alpha=0}^\infty
\lambda_\alpha}{\hb F(e^t)}).$$
By (\ref{prrod}) and the definition of $\overline{Z}_i$, we find
$$ 
\langle \phi_i, \overline{S}_{\mathbf{T}}(t,\hb) \rangle = 
e^{\lambda_i t/\hb} l\lambda_i
\overline{Z}_i(e^t, \hb,l).$$
Consider the correlator $S^*_\mathbf{T}(t,\hb,m+1)\in H_\mathbf{T}^*(\proj^m)
[[\hb^{-1},t,e^t]]$ defined
by (\ref{qertt}).
By equation (\ref{ertt}), 
the equality $\overline{Z}_i(e^t,\hb)=Z_i^*(e^t, \hb)$, 
 and property (\ref{vbvb}),
we conclude
$\overline{S}_{\mathbf{T}}(t,\hb) =S^*_\mathbf{T}(t,\hb,m+1)$.

After passing from equivariant to standard cohomology
$(\lambda_i=0)$ and setting $\hb=1$, we obtain the
Mirror result (case (iii) of Section 0). The series
$S_X^*, S_X \in H^*(\proj^m)[t][[e^t]]$ are determined by:
\begin{eqnarray*}
S_X^* & = & \frac{1}{(m+1)H}\ S_{\mathbf{T}}^*(t,\hb,m+1)|_{\lambda_i=0, 
\ \hb=1} 
\\
             & = & \sum_{i=0}^{m-1} I_i(t) H^i. \\
S_X & =& \frac{1}{(m+1)H} \ S_{\mathbf{T}}(T,\hb,m+1)|_{\lambda_i=0, \
\hb=1}.
\end{eqnarray*}
where $I_i(t) \in \Q[t][[e^t]]$ (see definitions (\ref{qertt}) and 
(\ref{deffff})).
The following equalities hold:
$$I_0(t)= F(e^t), \ \ I_1/I_0\ (t) = 
t+ \frac{(m+1)(G_{m+1}(e^t)-G_1(e^t))}{F(e^t)}.$$
We have shown $S_X$ is obtained
from $\sum_{i=0}^{m-1} I_i/I_0\ (t)$ by the
change of variables $T= I_1/I_0\ (t)$.
The proof of this explicit transformation
between $S_X^*$ and $S_X$ completes   
case (iii) of Section 0.

\subsection{The quintic 3-fold}
Let $X\subset \proj^4$ be a quintic 3-fold.
The expected dimension of the moduli space of
rational curves in $X$ is 0 for all degrees. 
The correlator
$S_X$ is easily evaluated in terms of the Gromov-Witten
invariants $N_d$ of $X$  
directly from the definitions.
Let $\mathcal{F}= 5T^3/6 + \sum_{d>0} N_d e^{dT}$.
After setting $\hb=1$, we obtain from (\ref{deffff}): 
$$S_X=\frac{1}{5H} \cdot  \sum_{d\geq 0} e^{(H+d)T} 
e_{2*}(\frac{c_{\text{top}}(E_d)}{1-\psi_2}).$$
It is necessary to calculate (for $d>0$):
\begin{eqnarray*}
e_{2*}(\frac{c_{\text{top}}(E_d)}{1-\psi_2}) &= &
e_{2*}(c_{\text{top}}(E_d) +  c_{\text{top}}(E_d)\psi_2 +
c_{\text{top}}(E_d)\psi^2_2) \\
& = & H^3 \cdot \langle \tau_0(1) \ \tau_1(H) \rangle^X_d + H^4
\cdot \langle \tau_0(1) \ \tau_2(1) \rangle^X_d \\
& = & dN_d H^3 -2 N_d H^4. 
\end{eqnarray*}
The expansion in the first line is truncated for
dimension reasons. The first term vanishes. Finally, the
string, dilaton, and divisor equations are applied to conclude the
last line. This integral calculation appears in [LLY] and [Ki].
An algebraic calculation now yields:
$$S_X =  1 +TH + \frac{1}{5}\frac{d\mathcal{F}}{dT}H^2 +
\Big(\frac{1}{5}T\frac{d\mathcal{F}}{dT} - \frac{2}{5}\mathcal{F}
\Big)H^3.$$
After
accounting for multiple covers by equation (\ref{manmul}),
$S_X$ exactly equals the ride side of (\ref{msy}).
Equation (\ref{ddeq}) follows from the quantum differential
equation obtained from the $*_X$ product (for which
$S_X$ is part of a fundamental solution). 
The proven correspondence (iii)  implies:
$$\mathcal{F}(T(t))= \frac{5}{2} \Big( \frac{I_1}{I_0}(t)
\frac{I_2}{I_0}(t) - \frac{I_3}{I_0}(t) \Big),$$ 
which is the standard form of the Mirror prediction for
quintic 3-folds.

\vspace{+10 pt}
\noindent
Department of Mathematics \\
University of Chicago \\
5734 S. University Ave. 60637 \\
rahul@math.uchicago.edu

\begin{thebibliography}{[COGP]}

\bibitem[AM]{am} P. Aspinwall and D. Morrison, {\em Topological
field theory and rational curves}, Comm. Math. Phys. {\bf 151}
(1993), 245-262.



\bibitem[Be]{be} A. Beauville, {\em Quantum cohomology of
complete intersections}, Mathematical Physics, Analysis, and
Geometry {\bf 168} (1995), 384-398.


\bibitem[B]{c} K. Behrend, {\em Gromov-Witten invariants
in algebraic geometry}, Invent. Math. {\bf 127} (1997), 601-617.

\bibitem[BF]{bf} K. Behrend and B. Fantechi, {\em The intrinsic normal
cone},  Invent. Math. {\bf 128} (1997), 45-88.


\bibitem[BM]{bm} K. Behrend and Yu. Manin, {\em Stacks of stable
maps and Gromov-Witten invariants}, Duke J. Math. {\bf 85} (1996)
no. 1, 1-60.

\bibitem[BDPP]{bdpp}  G. Bini, C. De Concini, M. Polito, and C. Procesi,
{\em Givental's work}, preprint 1998.


\bibitem[COGP]{cogp} P. Candelas, X.  de
la Ossa, P.  Green and L.  Parkes, {\em A pair of Calabi-Yau manifolds
as an exactly soluble superconformal field theory}, Nuclear Physics
{\bf B359} (1991), 21-74.

\bibitem[CK]{ck} D. Cox and S. Katz, {\em 
Mirror symmetry and algebraic geometry}, to appear in 
AMS Surveys and Monographs in Mathematics.

\bibitem[D]{d} B. Dubrovin, {\em The geometry of 2D topological
field theories}, in {\em Integrable systems and quantum
groups}, LNM {\bf 1620}, Springer-Verlag, 1996, 120-348.

\bibitem[ES1]{ees} G. Ellingsrud and S. Str\o mme, 
{\em The number of twisted cubics on the general quintic threefold},
Math. Scand. {\bf 76} (1995), no. 1, 5-34.

\bibitem[ES2]{eess} G. Ellingsrud and S. Str\o mme, 
{\em Bott's formula and enumerative geometry}, Jour. AMS {\bf 9}
(1996), no. 1, 175-193.


\bibitem[F]{f} W. Fulton, {\em Intersection theory}, Springer-Verlag,
1980.


\bibitem[FP]{fp}
 W. Fulton and R. Pandharipande,
{\em Notes on stable maps and quantum cohomology}, in 
{\em Proceedings of symposia in pure  mathematics: Algebraic geometry 
Santa Cruz 1995}, (J. Koll\'ar, R. Lazarsfeld, D. Morrison, eds.),
Volume 62, Part 2, 45-96.


\bibitem[G1]{g1} A. Givental, {\em Equivariant Gromov-Witten invariants},
Int. Math. Res. Notices {\bf 13} (1996), 613-663.

\bibitem[G2]{g2} A. Givental, {\em A mirror theorem for toric complete
intersections}, preprint 1997.


\bibitem[G3]{g3} A. Givental, {\em Elliptic Gromov-Witten invariants
and the generalized mirror conjecture}, preprint 1998.

\bibitem[GP]{gp} T. Graber and R. Pandharipande, {\em Localization
of virtual classes}, Invent. Math. (to appear).

\bibitem[GMP]{gmp} B. Greene, D. Morrison, and R. Plesser, {\em
Mirror manifolds in higher dimension}, Comm. Math. Phys.
Vol. 173 (1995), 559-598.

\bibitem[K]{k} S. Katz, {\em  On the finiteness of rational 
curves on quintic threefolds},  Compositio Math. {\bf 60} 
(1986), 151-162.

\bibitem[KJ1]{kj1} S. Kleiman and T. Johnsen, 
{\em Rational curves of degree at most $9$ on a general quintic
threefold}, Comm. Algebra {\bf 24} (1996), no. 8, 2721-2753.

\bibitem[KJ2]{kj2} S. Kleiman and T. Johnsen,
 {\em Toward Clemens' conjecture in degrees between 10 and 24}, 
Serdica Math
J. {\bf 23} (1997), 131-142.



\bibitem[Ki]{ki} B. Kim, {\em Quantum hyperplane section theorem
for homogeneous spaces}, preprint 1998.

\bibitem[Ko]{ko}  M. Kontsevich, {\em Enumeration of rational curves via 
torus actions}, in {\em The moduli space of curves},
(R. Dijkgraaf, C. Faber, and G. van der Geer, eds.), Birkhauser,
1995, 335-368.



\bibitem[KM]{km}      M. Kontsevich and Yu. Manin, {\em Gromov-Witten 
classes, quantum cohomology, and enumerative geometry}, Commun. Math.
Phys. {\bf 164} (1994), 525-562.

\bibitem[LLY]{lly} B. Lian, K. Liu, and S.-T. Yau, {\em Mirror principle
I}, Asian J. Math. Vol. 1, no. 4 (1997), 729-763.


\bibitem[LT]{lt} J. Li and G. Tian, {\em Virtual moduli cycles 
and Gromov-Witten invariants of algebraic varieties},
Jour. AMS {\bf{11}} (1998), no. 1, 119-174.



\bibitem[M]{m} Yu. Manin, {\em Generating functions in
algebraic geometry and sums over trees},
in {\em The moduli space of curves},
(R. Dijkgraaf, C. Faber, and G. van der Geer, eds.), Birkhauser,
1995, 401-417.


\bibitem[Mo]{mo} D. Morrison, {\em Mirror symmetry and rational
curves on quintic threefolds: A guide for mathematicians}, Jour. 
AMS {\bf 6} (1993), 223--247,


\bibitem[RT]{rt} Y. Ruan and G. Tian, {\em A mathematical
theory of quantum cohomology}, J. Diff. Geom. {\bf 42}
(1995), 259-367.

\bibitem[V]{v} I. Vainsencher,  {\em Enumeration of n-fold 
tangent hyperplanes to a surface},
J. Alg. Geometry,  {\bf 4} (1995), 503--526.


\bibitem[W1]{w1} E. Witten, {\em 
Mirror manifolds and topological field theory}, in:
{\em Essays on Mirror Manifolds}
(S.-T.\ Yau, ed.), International Press, Hong Kong 1992, 120-159.

\bibitem[W2]{w2} E. Witten, {\em Two dimensional gravity and intersection
theory on moduli space}, Survey in Diff. Geom. {\bf 1} (1991), 243-310.

\end{thebibliography}
\end{document}